%% file: sqnormsharp.tex
\documentclass[11pt,twoside]{article}
\newcommand{\ifims}[2]{#1} 
\newcommand{\ifAMS}[2]{#1}   
\newcommand{\ifau}[4]{#1}  
\newcommand{\ifbook}[2]{#1}   
\newcommand{\ifapp}[2]{#2}  
\newcommand{\ifadap}[2]{#1}  

\input myfrontmatter

    \renewcommand{\Section}[1]{\section{#1}}
    \renewcommand{\Subsection}[1]{\subsection{#1}}

\input mydef

\input statdef

\def\thetitle{Concentration of a high dimensional sub-gaussian vector}
\def\theruntitle{Concentration of a high dimensional sub-gaussian vector}

\def\theabstract{
This note describes the concentration property for a high dimensional sub-gaussian vector \( \Xv \).
In the Gaussian case, for any linear operator \( \QP \), it holds
\( \P\bigl( \| \QP \Xv \|^{2} - \tr (\BBH) > 2 \sqrt{\xx \tr(\BBH^{2})} + 2 \xx \| \BBH \| \bigr) \leq \ex^{-\xx} \) and
\( \P\bigl( \| \QP \Xv \|^{2} - \tr (\BBH) < - 2 \sqrt{\xx \tr(\BBH^{2})} \bigr) \leq \ex^{-\xx} \)
with \( \BBH = \QP \Var(\Xv) \QP^{\T} \); 
see \cite{laurentmassart2000}.
This implies concentration of the squared norm \( \| \QP \Xv \|^{2} \) around its expectation 
\( \E \| \QP \Xv \|^{2} = \tr (\BBH) \)
provided that \( \tr(\BBH^{2})/\| \BBH \|^{2} \) is sufficiently large.
An extension of this result to a non-gaussian case is a nontrivial task even under sub-gaussian behavior of \( \Xv \), especially if the entries of \( \Xv \) cannot be assumed independent 
and Hanson-Wright type bounds do not apply.
The results of this paper extend the Gaussian deviation bounds and support the concentration phenomenon 
for \( \| \QP \Xv \|^{2} \) using recent advances in Laplace approximation from \cite{SpLaplace2022}
and \cite{katsevich2023tight}.
The results are illustrated by the case when \( X \) is an i.i.d. sum.
}

\def\kwdp{62E15}
\def\kwds{62E10}

\def\thekeywords{exponential moments, quadratic forms, Laplace approximation}

\def\thankstitle{}

\aua
{Vladimir Spokoiny}
{Spokoiny, V. }
{
    Weierstrass Institute and HU Berlin,  
    \\
    Mohrenstr. 39, 10117 Berlin, Germany
    }
{spokoiny@wias-berlin.de}
{Weierstrass-Institute and Humboldt University Berlin}
{
  Financial support by the German Research Foundation (DFG) through the Collaborative Research Center 1294 ``Data assimilation'' is gratefully acknowledged.
}

\begin{document}
\thispagestyle{empty}
{
\title{\thetitle}
\theauthors

\maketitle
\begin{abstract}
{\footnotesize \theabstract}
\end{abstract}

\ifAMS
    {\par\noindent\emph{AMS Subject Classification:} Primary \kwdp. Secondary \kwds}
    {\par\noindent\emph{JEL codes}: \kwdp}

} 

\tableofcontents

\Section{Introduction}
\label{Sdevintr}
Let \( \xiv \sim \ND(0,\BBH) \) be a Gaussian vector in \( \R^{\dimp} \) with 
a covariance matrix \( \BBH \).
Also, let \( \tr(\BBH) \) denotes the trace of \( \BBH \) while \( \| \BBH \| \) is its operator norm.
The squared Euclidean norm \( \| \xiv \|^{2} \) fulfills \( \E \| \xiv \|^{2} = \tr (\BBH) \), 
\( \Var(\| \xiv \|^{2}) = \tr(\BBH^{2}) \), 
and 
for any \( \xx > 0 \)
\begin{EQ}[lcl]
	\P\Bigl( \| \xiv \|^{2} - \tr (\BBH) > 2 \sqrt{\xx \tr(\BBH^{2})} + 2 \xx \| \BBH \|  \Bigr)
	& \leq &
	\ex^{-\xx} ,
	\\
	\P\Bigl( \| \xiv \|^{2} - \tr (\BBH) < - 2 \sqrt{\xx \tr(\BBH^{2})} \Bigr)
	& \leq &
	\ex^{-\xx} ; 
\label{v7hfte535ghewjfyw}
\end{EQ}
see e.g. \cite{laurentmassart2000} for similar results or Theorem~\ref{TexpbLGA} in the Appendix.
If the value \( \dimH = \tr (\BBH) \) is much larger than \( \sqrt{\tr (\BBH^{2})} \)
then we observe the concentration property: the difference \( \| \xiv \|^{2} - \dimH \) is smaller (in probability) than \( \dimH \) itself.
The upper bound in \eqref{v7hfte535ghewjfyw} can be extended to a sub-gaussian case
using the same technique based on Markov's inequality; see e.g. \cite{HKZ2012} or Section~\ref{SsubGausstail}.

In the case of a vector \( \xiv \) with independent entries, 
the concentration property of the squared norm \( \| \xiv \|^{2} \) has been established 
in \cite{HaWr1971}; see also 
\cite{RuVe2013} and references therein for further developments.
In recent years, a number of new results have been obtained in this direction.
We refer to \cite{KlZi2019} and \cite{sambale23} for an extensive overview and advanced results on Hanson-Wright type concentration inequalities for a vector with independent sub-gaussian components and for new upper bounds established by the entropy method. 
\cite{AW2013} and \cite{GSS2021} demonstrated the power of log-Sobolev and Poincar\'e inequalities 
in obtaining upper deviation bounds for Gaussian and sub-gaussian tensors of higher order. 
Note, however, that independence of the entries of \( \xiv \) is a restrictive assumption
which enables rather strong tools like decoupling \cite{Gine2012} or entropic arguments \cite{KlZi2019}.
The main motivating example of this study comes from statistical applications when 
\( \xiv \) is a so called ``score vector''.
For such applications, independence of the entries of \( \xiv \) is not natural 
and is fulfilled only in very special cases like a sequence space model.
On the other side, for typical statistical setups with \( n \) independent observations, 
the score vector \( \xiv \) has an additive structure and it is asymptotically normal. 
However, these results require \( \dimp \) fixed and \( n \) large.
Modern statistical applications lead to situations with \( \dimp \) large and \( n \) small or
moderate. 
A crucial question arising in such situations is about validity of statistical procedures and applicability of related
guarantee bounds in terms of so called ``critical dimension'' meaning a kind of relation between 
\( n \) and \( \dimp \).
A motivating example below explains how all such questions can be reduced 
to extending \eqref{v7hfte535ghewjfyw} to the case of a general non-gaussian 
vector \( \xiv \) which is not covered by Hanson-Wright type results.
Rather it is related to a Gaussian approximation of the sort
\begin{EQA}
	\P\bigl( \| \xiv \| > \zq(\BBH,\xx) \bigr)
	& \approx &
	\P\bigl( \| \xivt \| > \zq(\BBH,\xx) \bigr)
	\leq 
	\ex^{-\xx} ,
\label{v0vkvv8gv7ey7fh5e3}
\end{EQA}
where \( \xivt \sim \ND(0,\BBH) \) and 
\( \zq^{2}(\BBH,\xx) = \tr(\BBH) + 2 \sqrt{\xx \tr(\BBH^{2})} + 2 \xx \| \BBH \| \).
A uniform bound like in \eqref{v7hfte535ghewjfyw} over all \( \xx \in \R_{+} \) is impossible 
because the large deviations of \( \| \QP \xiv \| \) involve the subgaussian constant. 
However, for many applications e.g. in statistics, biology, medicine, and econometrics, 
typical values of \( \alpha = \ex^{-\xx} \) belong to the range \( [0.01,0.05] \) 
which is covered by moderate values of \( \xx \).
The aim of this note is to provide some general conditions on \( \xiv \) ensuring nearly Gaussian quantiles of 
\( \| \xiv \|^{2} \). 
We establish a version of the upper bound in \eqref{v7hfte535ghewjfyw} using local smoothness of the moment generating
function \( \E \ex^{\langle \uv,\xiv \rangle} \) and the recent advances in Laplace approximation from 
\cite{SpLaplace2022} and \cite{katsevich2023tight}.
The lower bound is obtained by similar arguments applied to the characteristic function 
\( \E \ex^{\imi \langle \uv,\xiv \rangle} \).

\medskip

\textbf{A motivating example. Statistical estimation.}
Consider first a  linear model \( \Yv = \KSPsi \upsvs + \epsv \) in \( \R^{n} \), where 
\( \upsvs \in \R^{\dimp} \) is a target parameter, 
\( \epsv \) is a zero mean error vector in \( \R^{n} \), and 
\( \KSPsi \colon \R^{\dimp} \to \R^{n} \) is a given linear operator
such that \( \KSPsi^{\T} \KSPsi \) is invertible operator in \( \R^{\dimp} \). 
The prediction loss of the least squares solution \( \tilde{\upsv} = (\KSPsi^{\T} \KSPsi)^{-1} \KSPsi^{\T} \Yv \) satisfies
\begin{EQA}
	\| \KSPsi (\tilde{\upsv} - \upsvs) \|^{2}
	&=&
	\| \Pi \epsv \|^{2} \, ,
\label{8fijoiftgyefhjedkf}
\end{EQA}
where \( \Pi = \KSPsi (\KSPsi^{\T} \KSPsi)^{-1} \KSPsi^{\T} \).
This suggests the following construction of \emph{confidence sets}.
Given \( \xx \), fix \( \zq = \zq(\xx) \) by \( \P\bigl( \| \Pi \epsv \| > \zq \bigr) \leq \ex^{-\xx} \).
Then for \( \CA(\zq) \eqdef \{ \upsv \colon \| \KSPsi (\tilde{\upsv} - \upsv) \|^{2} \leq \zq^{2} \} \), it holds
\( \P(\upsvs \not\in \CA(\zq)) \leq \ex^{-\xx} \).
If \( \epsv \) is a Gaussian vector then the quantiles of \( \| \Pi \epsv \|^{2} \) can be
well approximated by \eqref{v7hfte535ghewjfyw} which only requires to know the covariance matrix \( \BBH = \Var(\Pi \epsv) \). 
The approximation result \eqref{v0vkvv8gv7ey7fh5e3} enables us to apply the same Gaussian quantiles 
even if \( \epsv \) is not Gaussian. 

More generally, consider a maximum likelihood estimator (MLE) \( \tilde{\upsv} \)
\begin{EQA}
	\tilde{\upsv}
	&=&
	\argmax_{\upsv \in \Ups} L(\upsv) ,
\label{9gnvyw37gyere8benxyw}
\end{EQA}
where \( \Ups \subset \R^{\dimp} \) is the parameter space and \( L(\upsv) \) is 
the log-likelihood function depending on the random data \( \Yv \in \R^{n} \).
Many other procedures like quasi MLE, minimum contract estimator, least absolute deviations (LAD), 
etc. can be considered similarly with \( L(\upsv) \) being negative loss function/contrast/fidelity.
This function may include some penalty for smoothness/complexity/regularity.
The target (background truth) \( \upsvs \) and the corresponding information operator \( \IF \) can be defined as
\begin{EQA}[c]
	\upsvs = \argmax_{\upsv} \E L(\upsv) ,
	\qquad
	\IF = \IF(\upsvs) = - \nabla^{2} \E L(\upsvs) \, .
\end{EQA}
\cite{Sp2022} offered a unified approach for studying the properties of \( \tilde{\upsv} \)
for a special class of so called stochastically linear smooth models for which the stochastic 
component \( \zeta(\upsv) = L(\upsv) - \E L(\upsv) \) is linear in \( \upsv \). 
This condition looks very restrictive but in many situations, it can be secured by an extension
of the parameter space; see e.g. \cite{Sp2023b}.
In particular, under this condition, with the \emph{standardized score vector} \( \xiv \eqdef \IF^{-1/2} \nabla \zeta \), it holds
on a set of dominating probability 
\begin{EQA}
	\| \IF^{1/2} (\tilde{\upsv} - \upsvs) - \xiv \|
	& \leq &
	\frac{3\dltwu_{3}}{4} \| \xiv \|^{2} \, ,
	\\
	\bigl| L(\tilde{\upsv}) - L(\upsvs) - \| \xiv \|^{2}/2 \bigr|
	& \leq &
	\dltwu_{3} \| \xiv \|^{3} \, ,
\label{odfodf0i30ioefkoeg}
\end{EQA}
where \( \dltwu_{3} \) is a small constant depending on local smoothness of \( \E L(\upsv) \) in a vicinity of \( \upsvs \);
see \cite{Sp2024}.
Therefore, if \( \dltwu_{3} \| \xiv \| \) is small, then
sharp deviation bounds for \( \| \xiv \| \) yield accurate concentration bounds and finite sample guarantees
for the errors \( \tilde{\upsv} - \upsvs \) and likelihood-based
confidence sets \( \CA(\zz) = \{ \upsv \colon L(\tilde{\upsv}) - L(\upsv) \leq \zz \} \).
These results provide a finite sample extension of the classical asymptotic results of parametric statistics;
see e.g. \cite{IH1981} or \cite{vdV00}.
In the case when the data \( \Yv \) are given as a sample of \( n \) independent observations,
the information matrix \( \IF \) scales with \( n \) while \( \dltwu_{3} \) is of order \( n^{-1/2} \).
Modern statistical tasks lead to the setup with a high dimensional vector \( \xiv \)
under limited sample size \( n \). 
%
The main result of the paper states a nonasymptotic version of \eqref{v7hfte535ghewjfyw} 
under condition of a large ``effective dimension'' \( \dimL = \tr(\BBH) \) and ``critical dimension'' condition
\( \dimL^{2} \ll n \).
%

The paper is organized as follows.
Section~\ref{SsubGausstail} provides a simple but rough upper bound under sub-gaussian condition on \( \Xv \).
The main results about concentration of \( \| \Xv \|^{2} \) are collected in Section~\ref{Sdevboundsharp}.
Section~\ref{Ssumiiddb} specifies the results to the case when \( \Xv \) is a normalized sum of independent random vectors,
and applicability of the presented bounds are discussed in Section~\ref{Srangesqnorm}.
The proofs of the main results from Section~\ref{Sdevboundsharp} rely on the results on
a high-dimensional Laplace approximation from Section~\ref{SlocalLaplace} and are deferred to Section~\ref{SqfnGproofs}.
%
Some useful technical facts about Gaussian quadratic forms are collected in Appendix~\ref{SmomentqfG} and Appendix~\ref{SdevboundGauss}.
Appendix~\ref{SDB3tens} provides some exponential bounds for a third-order Gaussian tensor
recently obtained in \cite{katsevich2023tight}.
These bounds are used in Section~\ref{SlocalLaplace}.

%
%

\input qf-nonGauss4

\input localLaplace

\input qf-nonGauss4proof

\newpage
\appendix

\input qf-Gauss
\input t3-Gauss4

\bibliography{exp_ts,listpubm-with-url}

\end{document}

%% file: myfrontmatter.tex
\usepackage{ifthen}
\usepackage[ampersand]{easylist}
\usepackage{amsmath,amssymb,amsthm}
\usepackage[bb=stix]{mathalpha}
\usepackage{mathtools}
\usepackage{natbib}
\usepackage{epsfig,graphicx}
\usepackage{comment}
\usepackage{color}
\usepackage{srcltx}
\usepackage[mathscr]{eucal}
\usepackage[math]{easyeqn}
\usepackage{etoolbox}
\usepackage{hyperref}
\hypersetup{
            colorlinks,
            linkcolor=hookersgreen,
            linktoc=blue,
            citecolor=ultramarine,
            urlcolor=black,
            filecolor=black
            }
\usepackage{nameref}
\usepackage{multirow}

\usepackage{accents}
\usepackage{mathrsfs}
\usepackage{bbm}
\usepackage{algorithm}
\usepackage{algpseudocode}

\ifims{
\textheight=22cm
\textwidth=14.8cm
\topmargin=0pt
\oddsidemargin=1.0cm
\evensidemargin=1.0cm
\linespread{1.3}

}{ 
}

\numberwithin{equation}{section}
\numberwithin{figure}{section}
\newcounter{example}[section]
\numberwithin{example}{section}
\newcounter{remark}[section]
\numberwithin{remark}{section}
\newtheorem{theorem}{Theorem}[section]
\newtheorem{proposition}[theorem]{Proposition}
\newtheorem{lemma}[theorem]{Lemma}
\newtheorem{corollary}[theorem]{Corollary}

\newtheorem{exmp}[example]{Example}
\newtheorem{rmrk}[remark]{Remark}
\newenvironment{example}{\begin{exmp}\rm}{\end{exmp}}
\newenvironment{remark}{\begin{rmrk}\rm}{\end{rmrk}}

\ifbook{
    
    \newcommand{\Section}[1]{\subsection{#1}}
    \newcommand{\Subsection}[1]{\subsubsection{#1}}

  }
  {
    
    \newcommand{\Section}[1]{\section{#1}}
    \newcommand{\Subsection}[1]{\subsection{#1}}

  }

%% file: mydef.tex
\renewcommand{\(}{$\,}
\renewcommand{\)}{\,$}

\def\nquad{\hspace{-1cm}}
\def\eqdef{\stackrel{\operatorname{def}}{=}}
\def\eqd{\stackrel{\operatorname{d}}{=}}

\DeclareMathAlphabet{\mathbbmsl}{U}{bbm}{bx}{sl}

\DeclareMathSymbol{\Alpha}{\mathalpha}{operators}{"41}
\DeclareMathSymbol{\Beta}{\mathalpha}{operators}{"42}
\DeclareMathSymbol{\Epsilon}{\mathalpha}{operators}{"45}
\DeclareMathSymbol{\Zeta}{\mathalpha}{operators}{"5A}
\DeclareMathSymbol{\Eta}{\mathalpha}{operators}{"48}
\DeclareMathSymbol{\Iota}{\mathalpha}{operators}{"49}
\DeclareMathSymbol{\Kappa}{\mathalpha}{operators}{"4B}
\DeclareMathSymbol{\Mu}{\mathalpha}{operators}{"4D}
\DeclareMathSymbol{\Nu}{\mathalpha}{operators}{"4E}
\DeclareMathSymbol{\Omicron}{\mathalpha}{operators}{"4F}
\DeclareMathSymbol{\Rho}{\mathalpha}{operators}{"50}
\DeclareMathSymbol{\Tau}{\mathalpha}{operators}{"54}
\DeclareMathSymbol{\Chi}{\mathalpha}{operators}{"58}
\DeclareMathSymbol{\omicron}{\mathord}{letters}{"6F}

\newcommand{\cc}[1]{\mathscr{#1}}
\newcommand{\bb}[1]{\boldsymbol{#1}}

\DeclareFontFamily{U}{mathx}{\hyphenchar\font45}
\DeclareFontShape{U}{mathx}{m}{n}{
<5><6><7><8><9><10>
<10.95><12><14.4><17.28><20.74><24.88>
mathx10
}{}
\DeclareSymbolFont{mathx}{U}{mathx}{m}{n}
\DeclareFontSubstitution{U}{mathx}{m}{n}
\DeclareMathAccent{\widebar}{0}{mathx}{"73}

\renewcommand{\bar}[1]{\widebar{#1}}
\renewcommand{\hat}[1]{\widehat{#1}}
\renewcommand{\tilde}[1]{\widetilde{#1}}

\makeatletter
\def\mathcenterto#1#2{\mathclap{\phantom{#1}\mathclap{#2}}\phantom{#1}}
\let\old@widetilde\widetilde
\def\widetildeto#1#2{\mathcenterto{#2}{\old@widetilde{\mathcenterto{#1}{#2\,}}}}
\let\old@widehat\widehat
\def\widehatto#1#2{\mathcenterto{#2}{\old@widehat{\mathcenterto{#1}{#2\,}}}}
\makeatother

\newcommand{\thankstitle}[1]{\ifthenelse{\equal{#1}{}}{}{\thanks{#1}}}
\newcommand{\thanksau}[1]{\ifthenelse{\equal{#1}{}}{}{\thanks{#1}}}

\newcommand{\aua}[6]
{\def\authora{#1}
\def\runauthora{#2}
\def\addressa{#3}
\def\emaila{#4}
\def\affiliationa{#5}
\def\thanksa{#6}}

\def\theauthors{
\ifau{ 
  \author{
    \authora
    \thanksau{\thanksa}
    \\[5.pt]
    \addressa \\
    \texttt{ \emaila}
  }
}
{  
  \author{
    \authora
    \thanksau{\thanksa}
    \\[5.pt]
    \addressa \\
    \texttt{ \emaila}
    \and
    \authorb
    \thanksau{\thanksb}
    \\[5.pt]
    \addressb \\
    \texttt{ \emailb}
  }
}
{   
  \author{
    \authora
    \thanksau{\thanksa}
    \\[5.pt]
    \addressa \\
    \texttt{ \emaila}
    \and
    \authorb
    \thanksau{\thanksb}
    \\[5.pt]
    \addressb \\
    \texttt{ \emailb}
    \and
    \authorc
    \thanksau{\thanksc}
    \\[5.pt]
    \addressc \\
    \texttt{ \emailc}
  }
} {   
  \author{
    \authora
    \thanksau{\thanksa}
    \\[5.pt]
    \addressa \\
    \texttt{ \emaila}
    \and
    \authorb
    \thanksau{\thanksb}
    \\[5.pt]
    \addressb \\
    \texttt{ \emailb}
    \and
    \authorc
    \thanksau{\thanksc}
    \\[5.pt]
    \addressc \\
    \texttt{ \emailc}
    \and
    \authord
    \thanksau{\thanksd}
    \\[5.pt]
    \addressd \\
    \texttt{ \emaild}
  }
}
}

\renewcommand{\Gamma}{\varGamma}
\renewcommand{\Pi}{\varPi}
\renewcommand{\Sigma}{\varSigma}
\renewcommand{\Delta}{\varDelta}
\renewcommand{\Lambda}{\varLambda}
\renewcommand{\Psi}{\varPsi}
\renewcommand{\Phi}{\varPhi}
\renewcommand{\Theta}{\varTheta}
\renewcommand{\Omega}{\varOmega}
\renewcommand{\Xi}{\varXi}
\renewcommand{\Upsilon}{\varUpsilon}

\def\argmax{\operatornamewithlimits{argmax}}

\def\av{\bb{a}}

\def\ev{\bb{e}}

\def\uv{\bb{u}}
\def\vv{\bb{v}}
\def\wv{\bb{w}}
\def\xv{\bb{x}}

\def\Av{\bb{A}}

\def\Uv{\bb{U}}

\def\Xv{\bb{X}}
\def\Yv{\bb{Y}}

\def\epsv{\bb{\varepsilon}}

\def\gammav{\bb{\gamma}}

\def\omegav{\bb{\omega}}

\def\xiv{\bb{\xi}}

\def\Psiv{\bb{\Psi}}

\def\sumi{\sum_{i=1}^{n}}

\usepackage{color}

\definecolor{blue(pigment)}{rgb}{0.2, 0.2, 0.6}
\definecolor{ultramarine}{rgb}{0.07, 0.04, 0.56}
\definecolor{darkspringgreen}{rgb}{0.09, 0.45, 0.27}
\definecolor{hookersgreen}{rgb}{0.0, 0.44, 0.0}
\definecolor{plum(traditional)}{rgb}{0.56, 0.27, 0.52}
\definecolor{purple(html/css)}{rgb}{0.5, 0.0, 0.5}
\definecolor{magenta(dye)}{rgb}{0.79, 0.08, 0.48}

%% file: statdef.tex
\def\GaussD{\gaussv_{\DPTG}}
\def\EX{\mathcal{E}}

\def\accu{\sigma}

\def\EUV{\E_{\, \UVL} \,}

\def\grad{\epsilon}

\def\VBH{S}
\def\VBHt{\tilde{\VBH}}

\def\leave{c}
\def\IIm{I^{\leave}}

\def\TG{\Gamma}
\def\TGD{\mathbbm{J}}

\def\trT{M}
\def\trTv{\bb{\trT}}
\def\trTvt{\tilde{\trTv}}

\def\perm{\pi}

\def\tensco{\tau}
\def\Tens{\mathcal{T}}
\def\Tenst{\tilde{\Tens}}
\def\TensG{\mathbb{T}}

\def\charf{\mathbb{f}}

\def\supH{\lambda}
\def\rexH{\epsilon}

\def\Egs{\E_{\gaussv}}
\def\Pgs{\P_{\gaussv}}

\def\hmax{\mathsf{c}}

\def\elll{\ell}

\def\lgd{f}

\def\lgdL{\elll}

\def\PfL{\P_{\lgd}}

\def\vH{\vA}

\def\WV{\mathcal{W}}

\def\Eta{\mathcal{H}}

\def\ic{i'}

\def\HVB{\mathbbmsl{V}} 

\def\dltw{\delta}

\def\dltwu{\tau}

\def\dltwhat{\hat{\dltw}}

\def\dlt{\delta}

\def\dltt{\tilde{\dlt}}

\def\vv{\bb{v}}

\def\II{\mathcal{I}}

\def\DFL{\mathbb{D}}
\def\DFLG{\DFL_{\GP}}

\def\R{\mathbbmsl{R}}
\def\E{\mathbbmsl{E}}

\def\P{\mathbbmsl{P}}

\def\kappa{\varkappa}

\def\diag{\operatorname{diag}}

\def\Fr{\operatorname{Fr}}

\def\ND{\mathcal{N}}

\def\Var{\operatorname{Var}}

\def\T{\top}
\def\tr{\operatorname{tr}}

\def\nsize{{n}}

\def\sumi{\sum_{i=1}^{\nsize}}

\def\ex{\mathrm{e}}

\def\Id{\mathbbmsl{I}}
\def\Ind{\operatorname{1}\hspace{-4.3pt}\operatorname{I}}

\def\alp{\alpha}

\def\BBH{B}

\def\cdens{\phi}

\def\CA{\mathcal{A}}

\def\CONST{\mathtt{C} \hspace{0.1em}}
\def\CONSTi{\mathtt{C}}

\def\CONSTgmb{\CONSTi_{\cdens}}
\def\CONSTchar{\CONSTi_{\charf}}

\def\const{\mathsf{c}}

\def\DP{D}

\def\DPH{\DP}
\def\DPTG{\mathsf{D}}

\def\DVL{\mathbb{D}}
\def\DVLG{\DVL_{\GP}}

\def\dimA{\mathbb{p}}

\def\dimL{\dimA}

\def\dimH{\dimA}

\def\dimp{p}

\def\dimQ{\dimA_{\QP}}

\def\dimq{q}

\def\Eta{\cc{H}}

\def\err{\diamondsuit}

\def\eps{\epsilon}			
\def\eps{\varepsilon}

\def\gaussv{\bb{\gauss}}
\def\gauss{\gamma}

\def\gaussG{\gaussv_{\GP}}

\def\gm{\mathtt{g}}

\def\gmb{\gm}

\def\gmn{\gm}
\def\GP{G}

\def\IF{\mathbbmsl{F}}

\def\IFL{{\mathbbmss{F}}}

\def\IFL{\mathbb{F}}

\def\imi{\mathtt{i}}

\def\jc{j'}

\def\KSPsi{\Psiv}

\def\Kappa{\cc{K}}

\def\kc{k'}

\def\muH{\mu}

\def\normG{\alpha}
\def\omegav{\bb{\phi}}

\def\QP{Q}

\def\rhogmn{\varrho}

\def\riskt{\cc{R}}

\def\rr{\mathtt{r}}

\def\rrL{\rr}

\def\rrTG{\rr}

\def\Tau{T}

\def\ups{\upsilon}
\def\upsv{\bb{\ups}}

\def\upsvs{\upsv^{*}}

\def\ups{\upsilon}
\def\upsv{\bb{\ups}}

\def\UV{\mathcal{U}}

\def\UVL{\UV}

\def\Ups{\varUpsilon}

\def\vA{\mathtt{v}}

\def\wv{\bb{w}}

\def\ww{w}

\def\xivt{\tilde{\xiv}}

\def\xvs{\xv^{*}}

\def\xx{\mathtt{x}}

\def\xxe{\xx_{\ex}}

\def\yy{\mathtt{y}}

\def\zq{z}

\def\zqe{\mathsf{z}}

\def\zz{\mathfrak{z}}

%% file: qf-nonGauss4.tex

\def\Xvt{\tilde{\Xv}}
\def\muHb{\bar{\muH}}
\def\gaussv{\gammav}
\def\rexH{\omega}
\def\CONSTgmb{\CONSTi_{\hspace{-1pt}X}}
\def\cdensX{\cdens_{\hspace{-1pt} X}}
\def\charfX{\charf_{\hspace{-1pt} X}}

\Section{Deviation bounds for sub-gaussian quadratic forms}
\label{Sprobabquad}
\label{SdevboundnonGauss}
This section collects some probability bounds for sub-gaussian quadratic forms.

\Subsection{A rough upper bound}
\label{SsubGausstail}

Let \( \xiv \) be a random vector in \( \R^{\dimp} \) with \( \E \xiv = 0 \).
We suppose that there exists a positive symmetric operator \( \HVB \) in \( \R^{\dimp} \) such that
\begin{EQA}
	\log \E \exp \bigl( \langle \uv, \HVB^{-1} \xiv \rangle \bigr)
	& \leq &
	\frac{\| \uv \|^{2}}{2} \, ,
	\qquad 
	\uv \in \R^{\dimp} .
\label{devboundinf}
\end{EQA}
In the Gaussian case, one can take \( \HVB^{2} = \Var(\xiv) \).
In general, \( \HVB^{2} \geq \Var(\xiv) \).
We consider a quadratic form \( \| \xiv \|^{2} \), where 
\( \xiv \) satisfies \eqref{devboundinf}.
We show that under \eqref{devboundinf}, the quadratic form \( \| \xiv \|^{2} \)
follows the same upper deviation bound 
\( \P\bigl( \| \xiv \|^{2} \geq \tr(\BBH) + 2 \sqrt{\xx \tr(\BBH^{2})} + 2 \xx \| \BBH \| \bigr) \leq \ex^{-\xx} \) 
with \( \BBH = \HVB^{2} \) as in the Gaussian case.

\begin{theorem}[\cite{HKZ2012}]
\label{Tdevboundinf}
Suppose \eqref{devboundinf}. 
Then for any \( \muH < 1/\| \BBH \| \)
\begin{EQA}
	\E \exp\Bigl( \frac{\muH}{2} \| \xiv \|^{2} \Bigr)
	& \leq &
	\exp \Bigl( \frac{\muH^{2} \tr(\BBH^{2})}{4 (1 - \| \BBH \| \muH)} + \frac{\muH \, \tr(\BBH)}{2} \Bigr)
\label{23iov96oklg7yf532tyfu}
\end{EQA}
and for any \( \xx > 0 \)
\begin{EQA}
	\P\bigl( \| \xiv \|^{2} > \tr(\BBH) + 2 \sqrt{\xx \tr(\BBH^{2})} + 2 \xx \| \BBH \| 
	\bigr)
	& \leq &
	\ex^{-\xx} .
\label{PxivbzzBBroBinf}
\end{EQA}
\end{theorem}

Statement \eqref{PxivbzzBBroBinf} looks identical to the upper bound in \eqref{v7hfte535ghewjfyw}, 
however, there is an essential difference:
\( \tr(\BBH) \) can be much larger than 
\( \E \| \xiv \|^{2} = \tr \Var(\xiv) \) if \( \HVB^{2} \) is much larger than \( \Var(\xiv) \).
The result from \eqref{PxivbzzBBroBinf} is not accurate enough
for supporting the concentration property that \( \| \xiv \|^{2} \) concentrates
around its expectation \( \E \| \xiv \|^{2} \).
The next section presents some sufficient conditions for obtaining sharp Gaussian-like deviation bounds. 

\Subsection{Concentration of the squared norm of a sub-gaussian vector} 
\label{Sdevboundsharp}
Let \( \Xv \) be a centered random vector in \( \R^{\dimp} \).
We study concentration property of the squared norm \( \| \QP \Xv \|^{2} \) for 
a linear mapping \( \QP \colon \R^{\dimp} \to \R^{\dimq} \).
The aim is to establish the results similar to 
\ifapp{\eqref{Pxiv2dimAvp12}}{\eqref{v7hfte535ghewjfyw}} 
with \( \BBH = \QP \Var(\Xv) \QP^{\T} \). 
Later we assume the following condition on the moment-generating function \( \E \ex^{ \langle \uv, \Xv \rangle } \).

\begin{description}
\item[\( \bb{(\cdensX)} \)\label{gmbref}]
\emph{A random vector \( \Xv \in \R^{\dimp} \) satisfies \( \E \Xv = 0 \), \( \Var(\Xv) \leq  \Id_{\dimp} \). 
The function \( \cdensX(\uv) \eqdef \log \E \ex^{ \langle \uv, \Xv \rangle } \) is finite and fulfills 
for some \( \CONSTgmb \)}
\begin{EQA}
	\cdensX(\uv)
	\eqdef
	\log \E \ex^{ \langle \uv, \Xv \rangle }
	& \leq &
	\frac{\CONSTgmb \| \uv \|^{2}}{2} \, ,
	\qquad
	\uv \in \R^{\dimp} \, .
\label{devboundinfgmb}
\end{EQA}
\end{description}

The condition \( \Var(\Xv) \leq  \Id_{\dimp} \) is only for convenience.
One can drop it by redefining \( \QP \) and \( \CONSTgmb \).
The constant \( \CONSTgmb \) can be quite large, it does not show up in the leading term of the obtained bound.
Also, we will only use \eqref{devboundinfgmb} for \( \| \uv \| \geq \gmn \) for some large \( \gmn \).
For \( \| \uv \| \leq \gmn \), we use smoothness properties of \( \cdensX(\uv) \).

The bounds in \ifapp{\eqref{Pxiv2dimAvp12}}{\eqref{v7hfte535ghewjfyw}} 
and in \eqref{PxivbzzBBroBinf} are uniform in the sense that they apply for all \( \xx \) and all \( \BBH \).
The results of this section are limited to a high dimensional situation
with \( \tr(\BBH^{2}) \gg \CONSTgmb \| \QP \QP^{\T} \| \) and apply only for 
\( \xx \ll \tr(\BBH^{2})/(\CONSTgmb \| \QP \QP^{\T} \|) \).
As compensation for this constraint, the bounds are surprisingly sharp.
In fact, they perfectly replicate bounds \ifapp{\eqref{Pxiv2dimAvp12}}{\eqref{v7hfte535ghewjfyw}} 
from the Gaussian case,
the upper and lower quantiles are exactly as in \ifapp{\eqref{Pxiv2dimAvp12}}{\eqref{v7hfte535ghewjfyw}} 
and the deviation probability 
is increased from  \( \ex^{-\xx} \) to \( (1 + \Delta_{\muH}) \ex^{-\xx} \) for a small value \( \Delta_{\muH} \).
For larger \( \xx \), one can still apply rough upper bound \eqref{PxivbzzBBroBinf} involving \( \CONSTgmb \). 

With \( \gaussv \) standard normal in \( \R^{\dimq} \), define
the \emph{effective trace} of \( \QP \) as
\begin{EQ}[rcl]
	\dimQ
	& \eqdef &
	\frac{\E \| \QP^{\T} \gaussv \|^{2}}{\| \QP \QP^{\T} \|}
	=
	\frac{\tr(\QP \QP^{\T})}{\| \QP \QP^{\T} \|} 
	=
	\frac{\| \QP \|_{\Fr}^{2}}{\| \QP \QP^{\T} \|} \, ,
\label{vrtgnfgih77jrfbegdhyd}
\end{EQ}
where \( \| \cdot \|_{\Fr} \) is the Frobenius norm.
%
For \( \wv \in \R^{\dimp} \), define a measure \( \P_{\wv} \) and the corresponding expectation \( \E_{\wv} \) such that for any r.v. \( \eta \)
\begin{EQA}
	\E_{\wv} \, \eta
	& \eqdef &
	\frac{\E (\eta \, \ex^{\langle \wv, \Xv \rangle})}{\E \ex^{\langle \wv, \Xv \rangle}} \, .	
\label{hcxuyjhgjvgui85ww3fg}
\end{EQA}
Also fix some \( \gmn > 0 \) and introduce 
\begin{EQA}[rcl]
	\dltwu_{3}
	& \eqdef &
	\sup_{\| \wv \| \leq \gmn} \, \sup_{\uv \in \R^{\dimp}}
	\frac{1}{\| \uv \|^{3}}
	\bigl| 
		\E_{\wv} \langle \uv, \Xv - \E_{\wv} \Xv \rangle^{3} 
	\bigr| \, ,
\label{7bvmt3g8rf62hjgkhgu3}
 	\\
	\dltwu_{4}
	& \eqdef &
	\sup_{\| \wv \| \leq \gmn} \, \sup_{\uv \in \R^{\dimp}} 
	\frac{1}{\| \uv \|^{4}}
	\bigl| 
		\E_{\wv} \langle \uv, \Xv - \E_{\wv} \Xv \rangle^{4} 
		- 3 \bigl\{ \E_{\wv} \langle \uv, \Xv - \E_{\wv} \Xv \rangle^{2} \bigr\}^{2} 
	\bigr| \, .
	\qquad
\label{7bvmt3g8rf62hjgkhgu}
\end{EQA}
The quantities \( \dltwu_{3} \) and \( \dltwu_{4} \) depend on the distribution of \( \Xv \) and \( \gmn \).
However, they are typically not only finite but also very small.
E.g. for \( \Xv \) Gaussian they just vanish.
If \( \Xv \) is a normalized sum of independent centred random vectors \( \xiv_{1},\ldots,\xiv_{n} \) then 
\( \dltwu_{3} \asymp n^{- 1/2} \) and \( \dltwu_{4} \asymp n^{- 1} \); see Section~\ref{Ssumiiddb}. 

First, we present an upper bound which extends 
\ifapp{\eqref{Pxiv2dimAvp12}}{\eqref{v7hfte535ghewjfyw}} to the non-Gaussian case. 

\begin{theorem}
\label{Tdevboundsharp}
Let \( \Xv \) satisfy \( \E \Xv = 0 \) and condition \nameref{gmbref}.
For any linear mapping \( \QP \colon \R^{\dimp} \to \R^{\dimq} \), define \( \BBH = \QP \Var(\Xv) \QP^{\T} \).
Let \( \gmn \) and \( \dltwu_{3} \) from \eqref{7bvmt3g8rf62hjgkhgu3} fulfill 
\( \gmn^{2} \geq 3 \dimQ \) and \( \gmn \, \dltwu_{3} \leq 2/3 \).
Then for any \( \xx > 0 \) with 
\( \sqrt{4 \xx } \leq \sqrt{\tr (\BBH^{2})}/(3 \CONSTgmb \| \QP \QP^{\T} \|^{2}) \), it holds
\begin{EQA}
	\P\bigl( \| \QP \Xv \|^{2} > \tr (\BBH) + 2 \sqrt{\xx \tr (\BBH^{2})} + 2\xx \| \BBH \| \bigr)
	& \leq &
	(1 + \Delta_{\muH}) \ex^{-\xx} \, ,
	\qquad
\label{PxivbzzBBroBinfsh}
\end{EQA}
where  \( \muH = \muH(\xx) \) is given by 
\( \muH^{-1} = \| \BBH \| + \sqrt{ \tr (\BBH^{2})/(4\xx)} \) and \( \Delta_{\muH} \) 
depends on \( \dltwu_{3} \), \( \dltwu_{4} \), \( \dimQ \) only and
will be given explicitly in the proof.
Moreover, \( \Delta_{\muH} \ll 1 \) under \( \tr (\BBH^{2}) \gg \| \QP \QP^{\T} \|^{2} \) and 
\( (\dltwu_{3}^{2} + \dltwu_{4}) \, \dimQ^{2} \ll 1 \).
\end{theorem}

\begin{remark}
The statement of Theorem~\ref{Tdevboundsharp} looks a bit technical, however, 
the main message is straightforward and useful:
for moderate \( \xx \)-values, the Gaussian upper quantiles \( \tr (\BBH) + 2 \sqrt{\xx \tr (\BBH^{2})} + 2\xx \| \BBH \| \)
ensure the nominal deviation probability \( \ex^{-\xx} \) even if \( \Xv \) is not Gaussian.
\end{remark}

For getting lower deviation bounds, 
in place of condition \nameref{gmbref} on the moment-generating function \( \E \exp\bigl( \langle \uv,\Xv \rangle \bigr) \),
we need a condition on the characteristic function \( \E \exp\bigl( \imi \langle \uv,\Xv \rangle \bigr) \).
Namely, we assume that it does not vanish and its logarithm is bounded on the ball \( \| \uv \| \leq \gmb \).

\begin{description}
\item[\( \bb{(\charfX)} \)\label{gmbiref}]
\emph{For some fixed \( \gmb \) and \( \CONSTchar \), the function 
\( \charfX(\uv) = \log \E \, \ex^{ \imi \langle \uv, \Xv \rangle } \) satisfies}
\begin{EQA}
	|\charfX(\uv)|
	=
	|\log \E \, \ex^{ \imi \langle \uv, \Xv \rangle }|
	& \leq &
	\CONSTchar \, ,
	\qquad
	\| \uv \| \leq \gmb \, .
\label{devboundinfgmbi}
\end{EQA}
\end{description}

Note that this condition can easily be ensured by replacing \( \Xv \) with \( \Xv + \alp \gaussv \) for any positive \( \alp \)
and \( \gaussv \sim \ND(0,\Id_{\dimp}) \).
The constant \( \CONSTchar \) is unimportant, it does not show up in our results.
It, however, enables us to define similarly to \eqref{7bvmt3g8rf62hjgkhgu}
\begin{EQA}[rcl]
	\dltwu_{4}
	& \eqdef &
	\sup_{\| \wv \| \leq \gmn} \,\, \sup_{\uv \in \R^{\dimp}}
	\frac{1}{\| \uv \|^{4}}
	\bigl| 
		\E_{\imi \wv} \langle \imi \uv, \Xv - \E_{\imi \wv} \Xv \rangle^{4} 
		- 3 \bigl\{ \E_{\imi \wv} \langle \imi \uv, \Xv - \E_{\imi \wv} \Xv \rangle^{2} \bigr\}^{2} 
	\bigr| \, .
	\qquad \quad
\label{7bvmt3g8rf62hjgkhgui}
\end{EQA}
The values \( \dltwu_{4} \) in \eqref{7bvmt3g8rf62hjgkhgu} and \eqref{7bvmt3g8rf62hjgkhgui} might be different,
however, we use the same notation without risk of confusion.

\begin{theorem}
\label{Texpquadroi}
Let \( \Xv \) satisfy \( \E \Xv = 0 \), \( \Var(\Xv) \leq  \Id_{\dimp} \).
Let also \( \QP \colon \R^{\dimp} \to \R^{\dimq} \) be a linear mapping,
\( \dimQ = \tr(\QP \QP^{\T}) \), \( \BBH = \QP \Var(\Xv) \QP^{\T} \). 
Assume \nameref{gmbiref} for some \( \gmn \) with \( \gmn^{2} \geq 3 \dimQ^{2} \).
Let also \( \dltwu_{3} \) be given by \eqref{7bvmt3g8rf62hjgkhgu3} and 
\( \gmn \, \dltwu_{3} \leq 2/3 \).
Then for any \( \xx \leq \tr(\BBH^{2})/4 \) 
\begin{EQA}
	\P\bigl( \| \QP \Xv \|^{2} < \tr (\BBH) - 2 \sqrt{\xx \tr(\BBH^{2})} \bigr)
	& \leq & 
	(2 + \err + \rho_{\muH}) \ex^{-\xx} \, ,
\label{PxivbzzBBroBinfshi}
\end{EQA}
where \( \muH \eqdef 2 \sqrt{\xx/\tr(\BBH^{2})} \) and
\begin{EQA}
	\rho_{\muH} 
	\eqdef
	\P\biggl( \| \QP^{\T} \gaussv \|^{2} \geq \frac{4 \muH^{-1} \dimQ^{2}}{\tr(\BBH^{2})} \biggr) 
	& \leq &
	\exp \Bigl( - \frac{4 \dimQ^{2}}{\tr(\BBH^{2})} \Bigr)
	\, .
	\qquad
\label{u8cnyhbnkjmjoyt8re3}
\end{EQA}
The value \( \err \) is given in the proof of Theorem~\ref{Texpquadro} and it is small
under \( \tr (\BBH^{2}) \gg \| \QP \QP^{\T} \| \) and \( (\dltwu_{3}^{2} + \dltwu_{4}) \, \dimQ \ll 1 \).
\end{theorem}

\Subsection{Sum of i.i.d. random vectors}
\label{Ssumiiddb}
Here we specify the obtained results to the case when \( \Xv = n^{-1/2} \sumi \xiv_{i} \) and 
\( \xiv_{i} \) are i.i.d. in \( \R^{\dimp} \) with \( \E \xiv_{i} = 0 \) and \( \Var(\xiv_{i}) = \Sigma \leq \Id_{\dimp} \).
In fact, only independence of the \( \xiv_{i} \)'s is used provided that all the moment conditions later on
are satisfied uniformly over \( i \leq n \). 
However, the formulation is slightly simplified in the i.i.d case.
Let some \( \QP \colon \R^{\dimp} \to \R^{\dimq} \) be fixed.
With \( \BBH = \QP \Sigma \QP^{\T} \), it holds \( \dimH = \E \| \QP \Xv \|^{2} = \tr (\BBH) \).
We study the concentration property for \( \| \QP \Xv \|^{2} \).
The goal is to apply Theorem~\ref{Tdevboundsharp} and Theorem~\ref{Texpquadroi} 
claiming that \( \| \QP \Xv \|^{2} - \dimH \) can be sandwiched between 
\( - 2 \sqrt{\xx \tr(\BBH^{2})} \) and \( 2 \sqrt{\xx \tr(\BBH^{2})} + 2\xx \| \BBH \| \) with probability at least \( 1 - 3 \ex^{-\xx} \).
The major required condition is sub-gaussian behavior of \( \xiv_{1} \).
The conditions are summarized here.

\begin{description}
\item[\( \bb{(\xiv_{1})} \)\label{gmb1ref}]
\emph{A random vector \( \xiv_{1} \in \R^{\dimp} \) satisfies \( \E \xiv_{1} = 0 \), \( \Var(\xiv_{1}) = \Sigma \leq  \Id_{\dimp} \). 
Also}
\begin{enumerate}
	\item \emph{The function \( \cdens_{\xiv}(\uv) \eqdef \log \E \ex^{ \langle \uv, \xiv_{1} \rangle } \) is finite and fulfills 
for some \( \CONSTgmb \)}
\begin{EQA}
	&&
	\cdens_{\xiv}(\uv)
	\eqdef
	\log \E \ex^{ \langle \uv, \xiv_{1} \rangle }
	\leq 
	\frac{\CONSTgmb \| \uv \|^{2}}{2} \, ,
	\qquad
	\uv \in \R^{\dimp} \, .
\label{devboundinfgmb1}
\end{EQA}
	\item
	\emph{For \( \rhogmn > 0 \) and some constants \( \hmax_{3} \) and \( \hmax_{4} \), it holds
	with \( \E_{\wv} \) from \eqref{hcxuyjhgjvgui85ww3fg}}
\begin{EQA}
	&&
	\sup_{\| \wv \| \leq \rhogmn} \,\, \sup_{\uv \in \R^{\dimp}}
	\frac{1}{\| \uv \|^{3}}
	\bigl| 
		\E_{\wv} \langle \uv, \xiv_{1} \rangle^{3} 
	\bigr| 
	\leq 
	\hmax_{3} \, ;
\label{d98k3efu7yvb67r4hfidke}
	\\
	&&
	\sup_{\| \wv \| \leq \rhogmn} \,\, \sup_{\uv \in \R^{\dimp}}
	\frac{1}{\| \uv \|^{4}}
	\bigl| 
		\E_{\wv} \langle \uv, \xiv_{1} - \E_{\wv} \xiv_{1} \rangle^{4} 
		- 3 \bigl\{ \E_{\wv} \langle \uv, \xiv_{1} - \E_{\wv} \xiv_{1} \rangle^{2} \bigr\}^{2} 
	\bigr| 
	\leq 
	\hmax_{4} \, .
\label{7bvmt3g8rf62hjgkhguiid}
\end{EQA}

	\item \emph{The function \( \log \E \, \ex^{ \imi \langle \uv, \xiv_{1} \rangle } \) is well defined and}
\begin{EQ}[c]
	\sup_{\| \wv \| \leq \rhogmn} \,\, \sup_{\uv \in \R^{\dimp}}
	\frac{1}{\| \uv \|^{4}}
	\bigl| 
		\E_{\imi \wv} \langle \imi \uv, \xiv_{1} - \E_{\imi \wv} \xiv_{1} \rangle^{4} 
		- 3 \bigl\{ \E_{\imi \wv} \langle \imi \uv, \xiv_{1} - \E_{\imi \wv} \xiv_{1} \rangle^{2} \bigr\}^{2} 
	\bigr| 
	\leq 
	\hmax_{4} \, .
\label{7bvmt3g8rf62hjgkhguiidi}
\end{EQ}
\end{enumerate}
\end{description}

We are now well prepared to state the result for the i.i.d. case. 

\begin{theorem}
\label{TnormXiid}
Let \( \Xv = n^{-1/2} \sumi \xiv_{i} \), where \( \xiv_{i} \) are i.i.d. in \( \R^{\dimp} \) satisfying
\( \E \xiv_{1} = 0 \) and \( \Var(\xiv_{1}) = \Sigma \leq \Id_{\dimp} \), and condition \nameref{gmb1ref}.
For a fixed \( \QP \), assume \( n \rhogmn^{2} \geq 3 \dimQ \) and \( n \gg \dimQ^{2} \).
Then with \( \BBH = \QP \Sigma \QP^{\T} \), it holds
\begin{EQA}[lcl]
	\P\Bigl( \| \QP \Xv \|^{2} - \tr (\BBH) > 2 \sqrt{\xx \tr (\BBH^{2})} + 2\xx \| \BBH \| \Bigr)
	& \leq &
	(1 + \Delta_{\muH}) \ex^{-\xx} \, ,
	\quad \text{ if } \quad
	\sqrt{4 \xx } \leq \frac{\sqrt{\tr (\BBH^{2})}}{3 \CONSTgmb \| \QP \QP^{\T} \|} \, ,
\label{hvb7ruehft6wjwiqws}
	\\
	\P\Bigl( \| \QP \Xv \|^{2} - \tr (\BBH) < - 2 \sqrt{\xx \tr (\BBH^{2})} \Bigr)
	& \leq &
	(2 + \Delta_{\muH}) \ex^{-\xx} \, ,
	\quad \text{ if } \quad
	\xx \leq \frac{\tr(\BBH^{2})}{4 \| \QP \QP^{\T} \|^{2}} ,
\label{gc82jwjv7w3f5wetjwy}
\end{EQA}
where \( \Delta_{\muH} \lesssim n^{-1} \dimQ^{2} \).
\end{theorem}

\Subsection{Range of applicability, critical dimension}
\label{Srangesqnorm}
This section discusses the range of applicability of the presented results, in particular, 
of the concentration property.
It was already mentioned earlier that concentration of the squared norm \( \| \QP \Xv \|^{2} \) 
is only possible in a high dimensional situation, even for \( \Xv \) Gaussian.
This condition can be written as \( \tr(\BBH^{2})/\| \QP \QP^{\T} \|^{2} \gg 1 \).
In our results, this condition is further detailed.
For instance, bound \eqref{PxivbzzBBroBinfsh} of Theorem~\ref{Tdevboundsharp} is only meaningful 
if \( \tr (\BBH^{2}) \gg \CONSTgmb^{2} \| \QP \QP^{\T} \|^{2} \).
This is the only place where the value \( \CONSTgmb \) shows up.

Another important quantity is the value \( \Delta_{\muH} \).
It should be small to make the presented results meaningful.
A sufficient condition for this property are
\( (\dltwu_{3}^{2} + \dltwu_{4}) \, \dimQ^{2} \ll 1 \).
For the i.i.d. case, this condition transforms into 
``critical dimension'' condition \( \dimQ^{2} \ll n \).
%
Recent results from \cite{katsevich2023tight} indicate that Laplace approximation could fail if 
\( \dimQ^{2} \ll n \) is not fulfilled even for a simple generalized linear model.
One can guess that a further relaxation of the ``critical dimension'' condition \( \dimQ^{2} \ll n \) is not possible and approximation 
\( \P\bigl( \| \QP \Xv \| > \zq(\BBH,\xx) \bigr) \approx \P\bigl( \| \QP \Xvt \| > \zq(\BBH,\xx) \bigr) \)
with a standard Gaussian vector \( \Xvt \) can fail if \( \dimQ^{2} \gg n \).


%% file: localLaplace.tex
\Section{Local Laplace approximation}
\label{SlocalLaplace}

This section presents some bounds on the error of local Laplace approximation.
Let \( \lgd(\xv) \) be a function in a high-dimensional Euclidean space \( \R^{\dimp} \) such that
\( \int \ex^{\lgd(\xv)} \, d\xv = \CONST < \infty \),
where the integral sign \( \int \) without limits means the integral over the whole space \( \R^{\dimp} \).
Then \( \lgd \) determines a distribution \( \PfL \) with the density
\( \CONST^{-1} \ex^{\lgd(\xv)} \).
Let \( \xvs \) be a point of maximum:
\begin{EQA}
	\lgd(\xvs)
	&=&
	\sup_{\uv \in \R^{\dimp}} \lgd(\xvs + \uv) .
\label{scdygw7ytd7wqqsquuqydtdtd}
\end{EQA}
We also assume that \( \lgd(\cdot) \) is at least three time differentiable. 
Introduce the negative Hessian \( \IFL = - \nabla^{2} \lgd(\xvs) \) and assume \( \IFL \) strictly positive definite.
Moreover, implicitly we assume that the negative Hessian \( \IFL = - \nabla^{2} \lgd(\xvs) \) is sufficiently large
in the sense that the Gaussian measure \( \ND(0,\IFL^{-1}) \) concentrates on a small local set \( \UVL \).
This allows to use a local Taylor expansion for 
\( \lgd(\xvs;\uv) \approx - \| \IFL^{1/2} \uv \|^{2}/2 \) in \( \uv \) on \( \UVL \).
For this local set \( \UVL \), we evaluate the quantity
\begin{EQA}
	\err
	& \eqdef &
	\biggl| \frac{\int_{\UVL} \ex^{\lgd(\xvs + \uv) - \lgd(\xvs)} \, d\uv - \int_{\UVL} \ex^{- \| \IFL^{1/2} \uv \|^{2}/2} \, d\uv} 
		 {\int \ex^{- \| \IFL^{1/2} \uv \|^{2}/2} d\uv} 
	\biggr| \, .
\label{fio9vkmfg763eu7g8gyhffr}
\end{EQA}
As \( \xvs = \argmax_{\xv} \lgd(\xv) \), it holds \( \nabla \lgd(\xvs) = 0 \) and 
\begin{EQA}
	\err
	& = &
	\biggl| \frac{\int_{\UVL} \ex^{\lgd(\xvs;\uv)} \, d\uv - \int_{\UVL} \ex^{- \| \IFL^{1/2} \uv \|^{2}/2} \, d\uv} 
		 {\int \ex^{- \| \IFL^{1/2} \uv \|^{2}/2} d\uv} 
	\biggr| ,
\label{IIgfifxutudufxutdu}
\end{EQA}
where \( \lgd(\xv;\uv) \) is the Bregman divergence 
\begin{EQA}
	\lgd(\xv;\uv)
	&=&
	\lgd(\xv + \uv) - \lgd(\xv) - \bigl\langle \nabla \lgd(\xv), \uv \bigr\rangle .
\label{fxufxpufxfpxu}
\end{EQA}
Our setup is motivated by Bayesian inference.
Assume that
\begin{EQA}
	\lgd(\xv)
	&=&
	\lgdL(\xv) - \| \GP (\xv - \xv_{0}) \|^{2} / 2
\label{jhdctrdfred4322edt7y}
\end{EQA}
for some \( \xv_{0} \) and a symmetric \( \dimp \)-matrix \( \GP^{2} \geq 0 \).
Here \( \lgdL(\cdot) \) stands for a log-likelihood function while the quadratic penalty 
\( \| \GP (\xv - \xv_{0}) \|^{2} / 2 \) corresponds to a Gaussian prior \( \ND(\xv_{0},\GP^{-2}) \).
Let also \( \DVL^{2} \eqdef - \nabla^{2} \lgdL(\xvs) > 0 \).
Then  
\begin{EQA}
	\IFL
	=
	- \nabla^{2} \lgd(\xvs)
	&=&
	- \nabla^{2} \lgdL(\xvs) + \GP^{2} 
	=
	\DVL^{2} + \GP^{2} .
\label{hydsf42wdsdtrdstdg}
\end{EQA} 
%
With decomposition \eqref{hydsf42wdsdtrdstdg} in mind, define 
\begin{EQA}[c]
	\DVLG^{2} = \IFL = \DVL^{2} + \GP^{2} \, ,
	\qquad
	\TGD^{2} 
	\eqdef 
	\DVLG^{-1} \DVL^{2} \DVLG^{-1} \, .
\label{dAdetrH02Hm2}
\end{EQA}
Assume that \( \lgd(\cdot) \) be a four times continuously differentiable function on \( \R^{\dimp} \).
Consider the remainder of the second and third-order Taylor approximation 
\begin{EQ}[rcl]
	\dltw_{3}(\uv)
	&=&
	\lgd(\xvs;\uv) - 
	\bigl\langle \nabla^{2} \lgd(\xvs) , \uv^{\otimes 2} \bigr\rangle/2 ,
	\\
	\dltw_{4}(\uv)
	&=&
	\lgd(\xvs;\uv) - 
	\bigl\langle \nabla^{2} \lgd(\xvs) , \uv^{\otimes 2} \bigr\rangle/2 
	- 
	\bigl\langle \nabla^{3} \lgd(\xvs), \uv^{\otimes 3} \bigr\rangle / 6 ,
\label{d4fuv1216303}
\end{EQ}
where \( \lgd(\xv;\uv) \) is given by \eqref{fxufxpufxfpxu}.
We will use the decomposition
\begin{EQA}
	\lgd(\xvs;\uv)
	&=&
	- \frac{1}{2} \| \DVLG^{-1} \uv \|^{2} + \dltw_{3}(\uv)
	=
	- \frac{1}{2} \| \DVLG^{-1} \uv \|^{2} + \Tens(\uv) + \dltw_{4}(\uv) ,
\label{yvjd7e3jgir63hfgkdd}
\end{EQA} 
where \( \Tens(\uv) = \langle \nabla^{3} \lgd(\xvs), \uv^{\otimes 3} \rangle/6 \) is the third order tensor
corresponding to the third derivative in the fourth order Taylor expansion for \( \lgd(\xvs;\uv) \).
For ease of notation, we skip dependence of \( \Tens \), \( \dltw_{3} \), and \( \dltw_{4} \) on \( \xvs \).

With some \( \rrL > 0 \), consider the local set \( \UVL \) of the form
\begin{EQA}
	\UVL 
	&=& 
	\bigl\{ \uv \colon \| \DVL \uv \| \leq \rrL \bigr\} .
\label{UvTDunm12spT}
\end{EQA}
Assume the following conditions.

\begin{description}
    \item[\label{l2l3sref} \( \bb{(\DVL_{3}^{*})} \)]
      \textit{For some \( \dltwu_{3} > 0 \), }
\begin{EQA}
	\sup_{\xv \colon \xv - \xvs \in \UVL} |\langle \nabla^{3} \lgd(\xv),\uv^{\otimes 3} \rangle|
	& \leq &
	\dltwu_{3} \, \| \DVL \uv \|^{3} ,
	\qquad
	\uv \in \R^{\dimp} \, .
\label{7cmvvc7e3hghjj856uied}
\end{EQA}
    \item[\label{l2l4ref} \( \bb{(\DVL_{4})} \)]
    \textit{For some \( \dltwu_{4} > 0 \), 
    }
\begin{EQA}
	|\dltw_{4}(\uv)|
	& \leq &
	\frac{\dltwu_{4}}{24} \| \DVL \uv \|^{4} ,
	\qquad
	\uv \in \UVL .
\label{vjb7hyer5ewre5fty4fgh}
\end{EQA}
\end{description}

Expansion \eqref{yvjd7e3jgir63hfgkdd} allows to represent 
the error \( \err \) from \eqref{IIgfifxutudufxutdu} as
\begin{EQA}
	\err
	=
	\frac{\int_{\UVL} \ex^{\lgd(\xvs;\uv)} \, d\uv - \int_{\UVL} \ex^{- \| \DVLG \uv \|^{2}/2} \, d\uv} 
		 {\int \ex^{- \| \DVLG \uv \|^{2}/2} d\uv} 
	&=&
	\EUV \Bigl[ \bigl\{ \exp \dltw_{3}(\gaussG) - 1 \bigr\} \Bigr] \, ,
\label{ed3le2d3pGd}
\end{EQA}
where \( \gaussG \sim \ND(0,\DVLG^{-2}) \) and \( \EUV \xi \) means \( \E \{ \xi \Ind(\gaussG \in \UVL) \} \).

\begin{proposition}
\label{PTensG42}
Assume \nameref{l2l3sref}, \nameref{l2l4ref}.
Define
\( \TGD^{2} = \DFLG^{-1} \DFL^{2} \, \DFLG^{-1} \),
\begin{EQA}[c]
	\grad = \dltwu_{3} \, \rr^{2} \| \TGD \|/2 ,
	\qquad
	\accu_{\GP}^{2} = \E \Tens^{2}(\gaussG) ,
	\qquad 
	\dltw_{4,\GP} = \EUV \dltw_{4}^{2}(\gaussG) .
\label{•}
\end{EQA}
Then 
\begin{EQA}[rcccl]
	\Bigl| \err - \frac{\accu_{\GP}^{2}}{2} \Bigr|
	& \leq &
	\accu_{\GP} \, \dltw_{4,\GP} 
	+ \frac{\dltw_{4,\GP}^{2}}{2} + \frac{5}{3} \grad^{3} \ex^{\grad^{2}} \, ,
	\quad
	\err
	& \leq &
	\frac{1}{2} (\accu_{\GP} + \dltw_{4,\GP})^{2} + \frac{5}{3} \grad^{3} \ex^{\grad^{2}} \, .
	\qquad
\label{u7yfjeoi8g8tj5udftrugh}
\end{EQA}
Moreover, 
\begin{EQA}[rcl]
\label{gdt6djiu97ie4425g7}
	\dltw_{4,\GP}
	& \leq &
	\frac{1}{24} \, \dltwu_{4} \, \bigl\{ \tr (\TGD^{2}) + 3 \| \TGD^{2} \| \bigr\}^{2} \, ,
	\\
	\accu_{\GP}
	& \leq &
	\sqrt{\frac{5}{12}} \,\, \dltwu_{3} \, \| \TGD \| \, \tr (\TGD^{2}) \, .
\label{gdt6djiu97ie4425g77}
\end{EQA}
\end{proposition}

\begin{proof}
%
%
Condition \nameref{l2l3sref} enables us to apply
Lemma~\ref{PtensGdlt} with \( X = \dltwhat_{3}(\gaussG) \) and \( k=3 \).
This yields 
\begin{EQA}
	\Bigl| 
		\EUV \Bigl\{ \Bigl(	\ex^{\dltwhat_{3}(\gaussG)} - 1 - \dltwhat_{3}(\gaussG) 
			- \frac{\dltwhat_{3}^{2}(\gaussG)}{2} \Bigr) g(\gaussG) 
		\Bigr\}
	\Bigr| 
	& \leq &
	\frac{5}{3} \grad^{3} \ex^{\grad^{2}} \, ,
\label{7nhig763hf5bv9edfuejn}
\end{EQA}
for any \( g(\cdot) \) with \( \sup_{\uv \in \UV} |g(\uv)| \leq 1 \).
Further, by \nameref{l2l4ref} and Lemma~\ref{Gaussmoments}
\begin{EQA}
	\EUV \dltw_{4}^{2}(\gaussG)
	& \leq &
	\frac{\dltwu_{4}^{2}}{24^{2}} \, \E \| \DFL \, \DFLG^{-1} \gaussv \|^{8}
	\leq 
	\frac{\dltwu_{4}^{2}}{24^{2}} \, \bigl\{ \tr(\TGD^{2}) + 3 \| \TGD^{2} \| \bigr\}^{4}
\label{yhdyvcyhwjg8it5j5uwfg}
\end{EQA}
and \eqref{gdt6djiu97ie4425g7} follows.
As \( \dltwhat_{3}(\gaussG) = \Tens(\gaussG) + \dltwhat_{4}(\gaussG) \), it holds
\begin{EQA}[rcl]
	\EUV |\dltwhat_{3}(\gaussG) - \Tens(\gaussG)| 
	&=& 
	\EUV |\dltwhat_{4}(\gaussG)|
	\leq 
	\sqrt{\EUV \dltwhat_{4}^{2}(\gaussG)}
	\leq 
	\sqrt{\EUV \dltw_{4}^{2}(\gaussG)} \, ,
\label{ycsy6ty25cgc53c7vyeetch}
\end{EQA}
and by \( \dltwhat_{4}(\gaussG) = \dltw_{4}(\gaussG) - \EUV \dltw_{4}(\gaussG) \)
\begin{EQA}[rcl]
	&& \nquad
	\EUV \bigl| \dltwhat_{3}^{2}(\gaussG) - \Tens^{2}(\gaussG) \bigr|
	\leq 
	2 \EUV \bigl| \dltwhat_{4}(\gaussG) \, \Tens(\gaussG) \bigr|
	+ \EUV \dltwhat_{4}^{2}(\gaussG)
	\\
	& \leq &
	2 \sqrt{\E \Tens^{2}(\gaussG) \, \EUV \dltwhat_{4}^{2}(\gaussG)}
	+ \EUV \dltwhat_{4}^{2}(\gaussG)
	\leq 
	2 \accu_{\GP} \, \dltw_{4,\GP} + \dltw_{4,\GP}^{2} \, .
\label{ycsy6ty25cgc53c7vyeetch2}
\end{EQA}
%
%
As \( \EUV \dltwhat_{3}(\gaussG) = 0 \), \eqref{7nhig763hf5bv9edfuejn} with \( g(\cdot) \equiv 1 \) 
and \eqref{ycsy6ty25cgc53c7vyeetch2} imply \eqref{u7yfjeoi8g8tj5udftrugh}.
The use of \eqref{uvjkb7nh3f6eyeyudt2jbi} from Lemma~\ref{LcolorGauss} with \( \tensco = \dltwu_{3}/6 \) yields
\begin{EQA}
	\EUV |\Tens(\gaussG)|
	& \leq &
	\sqrt{\EUV \Tens^{2}(\gaussG) }
	\leq 
	\frac{1}{6} \sqrt{15 \dltwu_{3}^{2} \, \| \TGD \|^{2} \, \tr^{2} (\TGD^{2})} 
\label{iujhdsufvhgfreu8478}
\end{EQA}
and \eqref{gdt6djiu97ie4425g77} follows as well. 
\end{proof}

%% file: qf-nonGauss4proof.tex
\ifadap{
\Section{Proofs of the main results}
\label{SqfnGproofs}
This section collects the proofs of the main results from Section~\ref{SdevboundnonGauss}.
}{}
\Subsection{Proof of Theorem~\ref{Tdevboundinf}}

Let \( \gaussv \) be standard Gaussian in \( \R^{\dimq} \) under \( \Egs \) conditionally on \( \xiv \).
For \( \muH \in (0,1) \), 
\begin{EQA}
	\E \exp\bigl( \muH \| \xiv \|^{2} / 2 \bigr)
	&=&
	\E \, \Egs \, \exp\bigl( \muH^{1/2} \langle \HVB \gammav, \HVB^{-1} \xiv \rangle \bigr) ,
\label{Egexmu12lHm1B12}
\end{EQA} 
Application of Fubini's theorem, \eqref{devboundinf}, and \eqref{m2v241m4mj1p} yields
\begin{EQA}
	\E \exp\Bigl( \frac{\muH}{2} \| \xiv \|^{2} \Bigr)
	& \leq &
	\Egs \, \exp\Bigl( \frac{\muH}{2} \| \HVB \gammav \|^{2} \Bigr)
	\leq 
	\exp \Bigl( \frac{\muH^{2} \tr (\BBH^{2})}{4 (1 - \muH \| \BBH \|)} + \frac{\muH \tr (\BBH)}{2} \Bigr) .
\label{wBmu2vA241mmutrBi}
\end{EQA}
Now the bound follows by Theorem~\ref{TexpbLGA} as in the Gaussian case.

\Subsection{Proof of Theorem~\ref{Tdevboundsharp}}
Normalizing by \( \| \QP \| \) reduces the statement to \( \| \QP \| = 1 \) and \( \dimQ = \tr(\QP \QP^{\T}) \) 
which will be supposed later.
This also implies \( \| \BBH \| = \| \QP \Var(\Xv) \QP^{\T} \| \leq 1 \).
The key step of the proof is the following statement.

\begin{proposition}
\label{Texpquadro}
Assume the conditions of Theorem~\ref{Tdevboundsharp} and \( \| \QP \| = 1 \).
If \( \muH > 0 \) satisfies 
\begin{EQA}[c]
	\CONSTgmb \, \muH \leq 1/3 
	\, ,
\label{uckdfyg6h6t43hfgy2}
\end{EQA}
then it holds 
\begin{EQA}
	\bigl| \E \exp ( \muH \| \QP \Xv \|^{2}/2) - \det (\Id_{\dimq} - \muH \BBH)^{-1/2} \bigr|
	& \leq &
	\Delta_{\muH} \det (\Id_{\dimq} - \muH \BBH)^{-1/2} \, 
	\qquad
\label{jcxu785t83w5ffr4ehjk}
\end{EQA}
for some constant \( \Delta_{\muH} \) such that \( \Delta_{\muH} \ll 1 \) under 
\( \dimQ \gg 1 \), \( (\dltwu_{3}^{2} + \dltwu_{4}) \, \dimQ^{2} \ll 1 \);
see the proof for a closed-form representation.
\end{proposition}

\begin{proof}
We use \eqref{Egexmu12lHm1B12} and Fubini theorem: with \( \Egs = \E_{\gauss \sim \ND(0,\Id_{\dimq})} \)
\begin{EQA}
	\E \exp\bigl( \muH \| \QP \Xv \|^{2} / 2 \bigr)
	&=&
	\E \, \Egs \, \exp\bigl( \muH^{1/2} \langle \QP^{\T} \gaussv, \Xv \rangle \bigr) 
	=
	\Egs \, \exp \cdensX(\muH^{1/2} \QP^{\T} \gaussv).
	\qquad
	\qquad
\label{Egexmu12lHm1B12X}
\end{EQA}
Further, redefine \( \gmn^{2} = 3 \dimQ \) and apply the decomposition
\begin{EQA}
	\Egs \, \exp \cdensX(\muH^{1/2} \QP^{\T} \gaussv)
	&=&
	\Egs \, \exp \cdensX(\muH^{1/2} \QP^{\T} \gaussv) \Ind(\| \muH^{1/2} \QP^{\T} \gaussv \| \leq \gmn)
	\\
	&& 
	+ \, \Egs \, \exp \cdensX(\muH^{1/2} \QP^{\T} \gaussv) \Ind(\| \muH^{1/2} \QP^{\T} \gaussv \| > \gmn) .
	\qquad
\label{pocvuw3jes45q4wtgdhjkj}
\end{EQA}
Each summand here will be bounded separately starting from the second one.
Define
\begin{EQA}[rcccl]
	\zz_{\muH} 
	& \eqdef &
	\frac{1}{4} \Bigl( \sqrt{\CONSTgmb^{-1} \muH^{-1} \gmn^{2}} - \sqrt{\dimQ} \Bigr)^{2} ,
	\qquad
	\rexH_{\muH} 
	& \eqdef &  
	\CONSTgmb \, \muH + \CONSTgmb \, \muH \sqrt{\dimQ / \zz_{\muH}} \, .
\label{jhdfuyewjwkifvgu4rtgy}
\end{EQA}
Then \eqref{uckdfyg6h6t43hfgy2} ensures that \( \zz_{\muH} \geq \bigl( \sqrt{9 \dimQ} - \sqrt{\dimQ} \bigr)^{2}/4 = \dimQ \) 
and \( \rexH_{\muH} \leq 2/3 \).
By \eqref{devboundinfgmb} and \eqref{llkknbononjm9hig4e} of Theorem~\ref{CTexpbLGA}, it holds under 
the condition \( \rexH_{\muH} \leq 2/3 \) 
\begin{EQA}
	&& \nquad
	\Egs \, \exp \cdensX(\muH^{1/2} \QP^{\T} \gaussv) \Ind(\| \muH^{1/2} \QP^{\T} \gaussv \| > \gmn)
	\\
	& \leq &
	\Egs \, \exp\bigl( \CONSTgmb \, \muH \| \QP^{\T} \gaussv \|^{2}/2 \bigr) 
	\Ind(\| \QP^{\T} \gaussv \|^{2} > \muH^{-1} \gmn^{2})
	\\
	& = &
	\exp \bigl( \CONSTgmb \, \muH \, \dimQ/2 \bigr) \,
	\Egs \, \exp\bigl( \CONSTgmb \, \muH (\| \QP^{\T} \gaussv \|^{2} - \dimQ)/2 \bigr) 
	\Ind(\| \QP^{\T} \gaussv \|^{2} > \muH^{-1} \gmn^{2})
	\\
	& \leq &
	\frac{1}{1 - \rexH_{\muH}} \, \exp\{ \CONSTgmb \, \muH \, \dimQ/2 - (1 - \rexH_{\muH}) \zz_{\muH} \} \, .
\label{hbcyedtcdfbghe653ddfh}
\end{EQA}
Note that \( \rexH_{\muH} \leq 2/3 \), \( \zz_{\muH} \geq \dimQ \), and \( \CONSTgmb \, \muH \leq 1/3 \) imply
\begin{EQA}
	\frac{1}{1 - \rexH_{\muH}} \, \exp\{ \CONSTgmb \, \muH \, \dimQ/2 - (1 - \rexH_{\muH}) \zz_{\muH} \}
	& \leq &
	3 \ex^{- \dimQ/6} \, .
\label{y7vcu8eiu36gvhrujeivq}
\end{EQA}
This inequality helps to bound the second term of \eqref{pocvuw3jes45q4wtgdhjkj} corresponding to the event 
\( \{ \| \muH^{1/2} \QP^{\T} \gaussv \| > \gmn \} \).
For the first term, we apply the results on Laplace approximation from Section~\ref{SlocalLaplace}.
First we check that \( \cdensX(\uv) \) satisfies conditions \nameref{l2l3sref} and \nameref{l2l4ref}:
\begin{EQA}
	|\langle \nabla^{3} \cdensX(\xv), \uv^{\otimes 3} \rangle|
	& \leq &
	\dltwu_{3} \| \uv \|^{3} \, , \quad
	\uv \in \R^{\dimp} ,
\label{jcu8dfuerg74jufuiuer3}
\end{EQA}
and
\begin{EQA}
	|\dltw_{4}(\uv)|
	& \eqdef &
	\Bigl| \cdensX(\uv) - \frac{1}{2} \langle \cdensX''(0) ,\uv^{\otimes 2} \rangle 
			- \frac{1}{6} \langle \cdensX^{(3)}(0) ,\uv^{\otimes 3} \rangle 
	\Bigr|
	\leq 
	\frac{\dltwu_{4}}{24} \| \uv \|^{4} \, , \quad
	\| \uv \| \leq \gmn \, .
	\qquad
\label{jcu8dfuerg74jufuiuer4}
\end{EQA}
Let us start with the univariate case.
Let \( X \) satisfy \( \E X = 0 \) and \( \E X^{2} \leq \sigma^{2} \).
Define for any \( t \in [0,\gmn] \)
a measure \( \P_{t} \) such that for any random variable \( \eta \)
\begin{EQA}
	\E_{t} \, \eta
	& \eqdef &
	\frac{\E (\eta \, \ex^{t X})}{\E \ex^{t X}} 
	\, .
\label{cu4e37gurreughierwo}
\end{EQA}
Consider \( \cdensX(t) \eqdef \log \E \ex^{t X} \) as a function of \( t \in [0,\gmn] \).
It is well defined and satisfies \( \cdensX(0) = \cdensX'(0) = 0 \),
\( \cdensX''(0) = \E X^{2} \leq \sigma^{2} \), and
\begin{EQA}
	\cdensX^{(3)}(t) 
	&=& \E_{t} (X - \E_{t} X)^{3} \, ,
	\\
	\cdensX^{(4)}(t)
	&=&
	\E_{t} (X - \E_{t} X)^{4} - 3 \bigl\{ \E_{t} (X - \E_{t} X)^{2} \bigr\}^{2} \, .	
\label{hyhboewhye7y6gvu45u8r}
\end{EQA}
Therefore, conditions \nameref{l2l3sref} and \nameref{l2l4ref} follow from 
\eqref{7bvmt3g8rf62hjgkhgu3} and \eqref{7bvmt3g8rf62hjgkhgu}.
The multivariate case can be reduced to the univariate one by fixing a direction \( \uv \in \R^{\dimp} \)
and considering the function \( \cdensX(t \uv) \) of \( t \).

Next, we apply Proposition~\ref{PTensG42} to evaluate the first term on the right-hand side of \eqref{pocvuw3jes45q4wtgdhjkj}.
Define \( \WV = \{ \wv \in \R^{\dimq} \colon \| \muH^{1/2} \QP^{\T} \wv \| \leq \gmn \} \).
Then with \( \gaussv \sim \ND(0,\Id_{\dimq}) \)
\begin{EQA}
	\Egs \, \exp \cdensX(\muH^{1/2} \QP^{\T} \gaussv) \Ind(\| \muH^{1/2} \QP^{\T} \gaussv \| \leq \gmn)
	& = &
	\CONSTi_{\dimq} \int_{\WV} \ex^{\lgd_{\muH}(\wv)} \, d\wv \, ,
\label{g6hnwyvhwweybgfewwel}
\end{EQA}
where \( \CONSTi_{\dimq} = (2\pi)^{-\dimq/2} \) and for \( \wv \in \R^{\dimq} \)
\begin{EQA}
	\lgd_{\muH}(\wv)
	&=&
	\cdensX(\muH^{1/2} \QP^{\T} \wv) - \| \wv \|^{2}/2
\label{vdsfbcttqbcfccsvcgarfcv}
\end{EQA}
so that \( \lgd_{\muH}(0) = 0 \), \( \nabla \lgd_{\muH}(0) = 0 \).
Also, define \( \DPH^{2} = \muH \QP \QP^{\T} \),
\begin{EQA}
\label{vthwctatqcuhqdgdsdsu3}
	\DPH_{\muH}^{2} 
	& \eqdef &
	- \nabla^{2} \lgd_{\muH}(0) 
	=
	- \muH \QP \Var(\Xv) \QP^{\T} + \Id_{\dimq}
	=
	\Id_{\dimq} - \muH \BBH ,
	\\
	\dimL_{\muH}
	& \eqdef &
	\tr \bigl\{ \DPH_{\muH}^{-1} (\muH \QP \QP^{\T}) \DPH_{\muH}^{-1} \bigr\} 
	=
	\muH \, \tr (\DPH_{\muH}^{-2} \, \QP \QP^{\T}),
	\\
	\normG_{\muH}
	& \eqdef &
	\| \DPH_{\muH}^{-1} (\muH \QP \QP^{\T}) \DPH_{\muH}^{-1} \| 
	=
	\muH \, \| \DPH_{\muH}^{-2} \, \QP \QP^{\T} \| \, .
\label{yyheu7futjthb884fhf6e}
\end{EQA}
By \eqref{uckdfyg6h6t43hfgy2} \( \muH \leq \frac{1}{3\CONSTgmb} \leq \frac{1}{3} \) and \( \| \BBH \| \leq \| \QP \| = 1 \) implies 
\( (1 - \muH) \Id_{\dimp} \leq \DPH_{\muH}^{2} \leq \Id_{\dimp} \) so that
\begin{EQA}
	\dimL_{\muH}
	& \leq &
	\frac{\muH}{1 - \muH} \, \dimQ  
	\leq 
	\frac{\dimQ}{2} \, ,
	\qquad
	\normG_{\muH}
	\leq  
	\frac{\muH}{1 - \muH} 
	\leq 
	\frac{1}{2}
	\, .
\label{ycv8w3ikdfbohyietwqr4cwq}
\end{EQA}
The function \( \lgd_{\muH}(\wv) \) inherits smoothness properties of \( \cdensX(\muH^{1/2} \QP^{\T} \wv) \).
In particular,  for any \( \wv \) with \( \| \muH^{1/2} \QP^{\T} \wv \| \leq \gmn \)
\begin{EQA}
	\bigl| \langle \nabla^{3} \lgd_{\muH}(\wv) , \uv^{\otimes 3} \rangle \bigr|
	& \leq &
	\dltwu_{3} \| \muH^{1/2} \QP^{\T} \uv \|^{3} 
	\, ,
\label{ucie939f7fgr5urfcews}
	\\
	\bigl| \langle \nabla^{4} \lgd_{\muH}(\wv) , \uv^{\otimes 4} \rangle \bigr|
	& \leq &
	\dltwu_{4} \| \muH^{1/2} \QP^{\T} \uv \|^{4} 
	\, .
\label{98cvkdy63ug7rtjhhh}
\end{EQA}
Now Proposition~\ref{PTensG42} applied to \( \lgd_{\muH}(\wv) \) yields
\begin{EQA}[c]
	\biggl| 
		\frac{\int_{\WV} \ex^{\lgd_{\muH}(\wv)} \, d\wv 
				- \int_{\WV} \ex^{- \| \DPH_{\muH} \wv \|^{2}/2} \, d\wv}
			 {\int \ex^{- \| \DPH_{\muH} \wv \|^{2}/2} \, d\wv}
	\biggr|
	\leq 
	\err \, .
\label{ghshyhcn2ynnsdjfgnso}
\end{EQA}
The quantity \( \err \) here is computed as follows.
Let \( \Tens(\uv) = \langle \nabla^{3} \lgd_{\muH}(0), \uv^{\otimes 3} \rangle \), 
\( \gaussv_{\muH} \sim \ND(0,\DPH_{\muH}^{-2}) \).
In view of \( \gmn^{2} = 3 \dimQ \) and \eqref{ycv8w3ikdfbohyietwqr4cwq}, it holds
\begin{EQA}
	\grad_{\muH}
	&=&
	\frac{\dltwu_{3} \, \gmn^{2} \sqrt{\normG_{\muH}}}{2} 
	\leq 
	\dltwu_{3} \, \dimQ \, ,
	\\
	\accu_{\muH}^{2} 
	&=& 
	\E \Tens^{2}(\gaussv_{\muH})
	\leq 
	\sqrt{5/12} \, \dltwu_{3} \, \dimL_{\muH}
	\leq 
	\frac{1}{3} \, \dltwu_{3} \, \dimQ \, ,
	\\
	\dltw_{4,\muH}
	& = &
	\EUV \dltw_{4}^{2}(\gaussv_{\muH})
	\leq 
	\frac{1}{24} \, \dltwu_{4} (\dimL_{\muH} + 3 \normG_{\muH})^{2} 
	\leq 
	\frac{1}{96} \, \dltwu_{4} (\dimQ + 3)^{2} \, .
\label{ydfje348bujtuy5r75t}
\end{EQA}
Then
\begin{EQA}[rcccl]
	\Bigl| \err - \frac{\accu_{\muH}^{2}}{2} \Bigr|
	& \leq &
	\accu_{\muH} \, \dltw_{4,\muH} 
	+ \frac{\dltw_{4,\muH}^{2}}{2} + \frac{5}{3} \grad_{\muH}^{3} \exp(\grad_{\muH}^{2}) \, ,
	\;\;
	\err
	& \leq &
	\frac{1}{2} (\accu_{\muH} + \dltw_{4,\muH})^{2} + \frac{5}{3} \grad_{\muH}^{3} \exp (\grad_{\muH}^{2}) \, .
	\qquad \quad
\label{vuf7fbfygtytewgetejj}
\end{EQA}
Furthermore, 
\begin{EQA}
	\rho_{\muH}
	& \eqdef &
	1 -	\frac{\int_{\WV} \ex^{- \| \DPH_{\muH} \wv \|^{2}/2} \, d\wv}{\int \ex^{- \| \DPH_{\muH} \wv \|^{2}/2} \, d\wv} 
	=
	\P\bigl( \| \sqrt{\muH} \QP^{\T} \DPH_{\muH}^{-1} \gaussv \| > \gmn \bigr)
	\, .
	\qquad
\label{yywnscthwsyd553t23e}
\end{EQA}
The use of \( \muH \leq 1/(3\CONSTgmb) \leq 1/3 \) and \( (1 - \muH) \muH^{-1} \gmn^{2} \geq 6 \CONSTgmb \, \dimQ \geq 6 \dimQ \) 
yields
\begin{EQA}
	\rho_{\muH}
	& \leq &
	\P\bigl( \| \QP^{\T} \gaussv \|^{2} > \frac{(1 - \muH) \gmn^{2}}{\muH} \bigr) 
	\leq 
	\P\bigl( \| \QP^{\T} \gaussv \|^{2} > 6 \dimQ \bigr)
	\leq 
	\ex^{- \dimQ/2} \, ;
\label{8cfkw9vire78vjweiowsck7}
\end{EQA}
see \eqref{3emzmsp22z2} of Theorem~\ref{CTexpbLGA} with \( \yy = 5 \dimQ \).
By \eqref{ghshyhcn2ynnsdjfgnso} and \eqref{yywnscthwsyd553t23e}
\begin{EQA}
	\biggl| \frac{\int_{\WV} \ex^{\lgd_{\muH}(\wv)} \, d\wv}{\int \ex^{- \| \DPH_{\muH} \wv \|^{2}/2} \, d\wv} - 1 \biggr|
	& \leq &
	\err + \rho_{\muH} \, .
\label{obmfurfwubyehjeuy}
\end{EQA}
It remains to be noted that
\begin{EQA}
	\CONSTi_{\dimq} \int \ex^{- \| \DPH_{\muH} \wv \|^{2}/2} \, d\wv
	&=&
	\frac{1}{\det \DPH_{\muH}}
	=
	\det(\Id_{\dimq} - \muH \BBH)^{-1/2} 
	\leq 
	1
\label{ivmeyh8gu467gsbewrde}
\end{EQA}
and \eqref{jcxu785t83w5ffr4ehjk} follows from \eqref{hbcyedtcdfbghe653ddfh} and \eqref{obmfurfwubyehjeuy} with
\begin{EQA}
	\Delta_{\muH}
	& \leq &
	\err + \rho_{\muH} + 
	\frac{1}{1 - \rexH_{\muH}} \, \exp\{ \CONSTgmb \, \muH \, \dimQ/2 - (1 - \rexH_{\muH}) \zz_{\muH} \} \, .
	\qquad
\label{jcxu785t83w5ffr4ehjkde}
\end{EQA}
Moreover,  for \( \dimQ \) large, \( \rho_{\muH} \) is small by \eqref{8cfkw9vire78vjweiowsck7},
the exp-term can be bounded by \eqref{y7vcu8eiu36gvhrujeivq},  
while \( \err \) is small provided that \( (\dltwu_{3}^{2} + \dltwu_{4}) \, \dimQ^{2} \) is small.
\end{proof}

Now we can finalize the proof of Theorem~\ref{Tdevboundsharp}.
Upper deviation bounds for \( \| \QP \Xv \|^{2} \) can be derived as in the Gaussian case 
by applying \eqref{jcxu785t83w5ffr4ehjk} with a proper choice of \( \muH \).
Let \( \xx \) satisfy \( \sqrt{4 \xx } \leq \sqrt{\tr (\BBH^{2})}/(3 \CONSTgmb \| \BBH \|) \).
We check \eqref{uckdfyg6h6t43hfgy2} for \( \muH = \muH(\xx) \).
Indeed, the definition \( \muH^{-1} = \| \BBH \| + \sqrt{ \tr (\BBH^{2})/(4\xx)} \) 
implies \( \muH \leq \sqrt{4 \xx / \tr (\BBH^{2})} \). 
Therefore, \( \sqrt{4 \xx } \leq \sqrt{\tr (\BBH^{2})}/(3 \CONSTgmb) \) yields
\( \muH \leq 1/(3 \CONSTgmb) \) and \eqref{uckdfyg6h6t43hfgy2} is fulfilled for \( \gmn^{2} = 3 \dimQ \).
The bound \eqref{PxivbzzBBroBinfsh} follows from \eqref{jcxu785t83w5ffr4ehjk} as in the Gaussian case of Theorem~\ref{TexpbLGA}.

\Subsection{Proof of Theorem~\ref{Texpquadroi}}

This result is based on an approximation 
\( \E \ex^{ - \muH \| \QP \Xv \|^{2}/2} \approx \det (\Id_{\dimq} + \muH \BBH)^{-1/2} \).
This requires an analog of \eqref{jcxu785t83w5ffr4ehjk} for \( \muH \) negative.
With \( \imi = \sqrt{-1} \), the use of \eqref{Egexmu12lHm1B12X} yields
\begin{EQA}
	\E \, \ex^{ - \muH \| \QP \Xv \|^{2} / 2} 
	&=&
	\E \, \Egs \, \ex^{ \imi \sqrt{\muH} \langle \QP^{\T} \gaussv, \Xv \rangle} 
	=
	\Egs \E \, \ex^{ \imi \sqrt{\muH} \langle \QP^{\T} \gaussv,\Xv \rangle} \, .
\label{Egexmu12lHm1B12Xch}
\end{EQA}

\begin{proposition}
\label{Pexpquadroi}
Assume the conditions of Theorem~\ref{Texpquadroi}.
For any \( \muH \in (0,1) \), it holds with \( \BBH = \QP \Var(\Xv) \QP^{\T} \)
\begin{EQ}[rcl]
	\bigl| \E \ex^{ - \muH \| \QP \Xv \|^{2}/2} - \det (\Id_{\dimq} + \muH \BBH)^{-1/2} \bigr|
	& \leq &
	(\err + \rho_{\muH}) \det (\Id_{\dimq} + \muH \BBH)^{-1/2} + \rho_{\muH} \, ,
	\\
	\rho_{\muH} 
	\leq  
	\Pgs\bigl( \| \QP^{\T} \gaussv \|^{2} \geq 4 \muH^{-1} \dimQ \bigr)
	& \leq &
	\exp \Bigl\{ - \frac{\dimQ}{4} (2\muH^{-1/2} - 1)^{2} \Bigr\} \, .
\label{jcxu785t83w5ffr4ehjki}
\end{EQ}
\end{proposition}

\begin{proof}
We follow the proof of Proposition~\ref{Texpquadro} replacing everywhere \( \cdensX(\uv) \) with \( \charfX(\uv) \).
In particular, we start with representation \eqref{Egexmu12lHm1B12Xch} and apply with \( \gmn^{2} = 3 \dimQ \)
\begin{EQA}
	&& \nquad
	\E \, \ex^{ - \muH \| \QP \Xv \|^{2} / 2}
	=
	\Egs \, \ex^{\charfX(\sqrt{\muH} \, \QP^{\T} \gaussv)}
	\\
	&=&
	\Egs \, \ex^{\charfX(\sqrt{\muH} \, \QP^{\T} \gaussv)} \Ind(\| \sqrt{\muH} \, \QP^{\T} \gaussv \| \leq \gmn)
	+ \Egs \, \ex^{\charfX(\sqrt{\muH} \, \QP^{\T} \gaussv)} \Ind(\| \sqrt{\muH} \, \QP^{\T} \gaussv \| > \gmn) .
\label{pocvuw3jes45q4wtgdhjkji}
\end{EQA}
It holds 
\begin{EQA}
	\charfX(0)
	&=&
	0,
	\quad
	\nabla \charfX(0)
	=
	0,
	\quad
	- \nabla^{2} \charfX(0) = \Var(\Xv) \leq \Id_{\dimp} \, .
\label{ycywhwf5e4g3eerter}
\end{EQA}
Moreover, smoothness conditions \eqref{jcu8dfuerg74jufuiuer3}, \eqref{jcu8dfuerg74jufuiuer4} 
are automatically fulfilled for \( \charfX(\uv) \)
with the same \( \dltwu_{3} \) and \( \dltwu_{4} \) from \eqref{7bvmt3g8rf62hjgkhgui}.
The most important observation for the proof is that the bound \eqref{obmfurfwubyehjeuy} continues to apply for 
\begin{EQA}
	\lgd_{\muH}(\wv)
	&=&
	\charfX(\sqrt{\muH} \, \QP^{\T} \wv) - \| \wv \|^{2}/2,
\label{ujc563hfg8h75urdwyehfcn}
\end{EQA}
with \( \err \) from \eqref{vuf7fbfygtytewgetejj} and 
\begin{EQA}[rcccl]
	\DPH_{\muH}^{2} 
	& \eqdef &
	- \nabla^{2} \lgd_{\muH}(0) 
	&=&
	\muH \QP \Var(\Xv) \QP^{\T} + \Id_{\dimq}
	=
	\Id_{\dimq} + \muH \BBH ,
	\\
	\dimL_{\muH}
	& \eqdef &
	\tr \bigl\{ \DPH_{\muH}^{-2} (\muH \QP \QP^{\T}) \bigr\}
	&\leq &
	\frac{\muH}{1 + \muH} \tr (\QP \QP^{\T}) 
	\leq 
	\muH \, \dimQ \, ,
	\\
	\normG_{\muH}
	& \eqdef &
	\| \DPH_{\muH}^{-1} (\muH \QP \QP^{\T}) \DPH_{\muH}^{-1} \| 
	&\leq &
	\frac{\muH}{1 + \muH} 
	\leq 
	\muH \, ,
\label{vthwctatqcuhqdgdsdsu3i}
\end{EQA}
and \( \rho_{\muH} \leq \P\bigl( \| \QP^{\T} \gaussv \|^{2} \geq \muH^{-1} \gmn^{2} \bigr) \); cf. \eqref{yywnscthwsyd553t23e}.
This yields
\begin{EQA}
	\biggl| \Egs \, \ex^{\charfX(\sqrt{\muH} \, \QP^{\T} \gaussv)} \Ind(\| \sqrt{\muH} \, \QP^{\T} \gaussv \| \leq \gmn) 
		- \frac{1}{\det (\Id_{\dimq} + \muH \BBH)^{1/2}} 
	\biggr|
	& \leq &
	\frac{\err + \rho_{\muH}}{\det (\Id_{\dimq} + \muH \BBH)^{1/2}} \, .
\label{0vkr7vhfh4twhswtfhry4y7}
\end{EQA}
Finally we use \( |\ex^{\charfX(\uv)}| \leq 1 \) and thus,
\begin{EQA}
	\bigl| \Egs \, \ex^{\charfX(\sqrt{\muH} \, \QP^{\T} \gaussv)} \Ind(\| \sqrt{\muH} \, \QP^{\T} \gaussv \| > \gmn) \bigr|
	& \leq &
	\P\bigl( \| \sqrt{\muH} \, \QP^{\T} \gaussv \| > \gmn \bigr) 
	=
	\Pgs\bigl( \| \QP^{\T} \gaussv \|^{2} \geq 4 \muH^{-1} \dimQ \bigr)
\label{hdveyhwedvuwuwejdv1334}
\end{EQA}
and \eqref{jcxu785t83w5ffr4ehjki} follows in view of \eqref{3emzmsp22z2} of Theorem~\ref{CTexpbLGA}.
\end{proof}

The proof of Theorem~\ref{Texpquadroi} is similar to the case of Theorem~\ref{Tdevboundsharp}.
By the exponential Chebyshev inequality and \eqref{jcxu785t83w5ffr4ehjki}, it holds with \( \vH^{2} = \tr(\BBH^{2}) \)
\begin{EQA}
	&& \nquad
	\P\bigl( \tr (\BBH) - \| \QP \Xv \|^{2} > 2 \vH \sqrt{\xx} \bigr)
	\leq
	\exp( - \muH \, \vH \sqrt{\xx}) \E \exp \bigl\{ \muH \tr (\BBH)/2 - \muH \| \QP \Xv \|^{2} /2 \bigr\}
	\\
	& \leq &
	\exp(\muH \tr (\BBH)/2 - \muH \, \vH \sqrt{\xx}) 
	\bigl\{ (1 + \err + \rho_{\muH}) \det (\Id_{\dimq} + \muH \BBH)^{-1/2} + \rho_{\muH} \bigr\} .
\label{cuyfr7ygbyuby6r4yfwete}
\end{EQA}
In view of \( x - \log(1+x) \leq x^{2}/2 \) and \( \muH = 2 \vH^{-1} \sqrt{\xx} \), 
as in the proof of Lemma~\ref{Lqfexpmom}
\begin{EQA}
	- \muH \, \vH \sqrt{\xx} + \muH \tr (\BBH)/2 + \log \det(\Id_{\dimq} + \muH \BBH)^{-1/2}
	& \leq &
	- \muH \, \vH \sqrt{\xx} + \muH^{2} \vH^{2}/4
	=
	- \xx .
\label{d7ywh3dfuyt523fhfyes}
\end{EQA}
Also \( \muH \tr (\BBH)/2 - \muH \, \vH \sqrt{\xx} = \vH^{-1} \tr (\BBH) \, \sqrt{\xx} - 2 \xx \leq \vH^{-1} \dimQ \, \sqrt{\xx} - 2 \xx \).
The bound on \( \rho_{\muH} \) in \eqref{u8cnyhbnkjmjoyt8re3} follows from \eqref{3emzmsp22z2} of Theorem~\ref{CTexpbLGA}
in view of \( \dimQ \geq \vH^{2} \) and hence, \( \dimQ \leq \dimQ^{2}/\vH^{2} \).
Finally, observe that
\begin{EQA}
	\rho_{\muH} \, \exp\bigl( \vH^{-1} \dimQ \, \sqrt{\xx} - 2 \xx \bigr)
	& \leq &
	\exp\Bigl( - \frac{\dimQ^{2}}{4\vH^{2}} + \frac{\dimQ \, \sqrt{\xx}}{\vH} - 2 \xx \Bigr)
	\\
	& \leq &
	\exp\Bigl\{ - \Bigl( \frac{\dimQ}{2 \vH} - \sqrt{\xx} \Bigr)^{2} - \xx \Bigr\}
	\leq 
	\ex^{-\xx} \, 
\label{juf7wjejfc6erdyrjvgbiw3}
\end{EQA}
and \eqref{PxivbzzBBroBinfshi} follows as well.

\Subsection{Proof of Theorem~\ref{TnormXiid}}
The definition  and i.i.d structure of the \( \xiv_{i} \)'s yield 
\( \E \langle \uv,\Xv \rangle^{2} =	\E \langle \uv,\xiv_{1} \rangle^{2} \) and 
\begin{EQA}
	\cdensX(\uv)
	&=&
	\log \E \ex^{\langle \Xv,\uv \rangle}
	=
	n \cdens_{\xiv}(n^{-1/2}\uv) 
\label{ufiu3jgibnhi5keuvytr}
\end{EQA}
for any \( \uv \in \R^{\dimp} \),
where \( \cdens_{\xiv}(\uv) \eqdef \log \E \ex^{\langle \xiv_{1},\uv \rangle} \). 
For the derivatives \( \cdensX^{(k)}(\uv) \) of \( \cdensX(\uv) \), this yields
\begin{EQA}
	\cdensX^{(k)}(\uv)
	&=&
	n^{1 - k/2} \cdens_{\xiv}^{(k)}(n^{-1/2}\uv) .
\label{ufiu3jgibnhi5keuvytr}
\end{EQA}
This enables us to derive \eqref{7bvmt3g8rf62hjgkhgu3} and \eqref{7bvmt3g8rf62hjgkhgu} from \nameref{gmb1ref}
for any \( \gmn \) with \( \gmn/\sqrt{n} \leq \rhogmn \) and
\begin{EQA}
	\dltwu_{3}
	& = &
	n^{-1/2} \hmax_{3} \, ,
	\qquad
	\dltwu_{4}
	= 
	n^{-1} \hmax_{4} \, .
\label{gv8iejkwidwweewfg}
\end{EQA}
Moreover, the quantity \( \err \) from \eqref{vuf7fbfygtytewgetejj} satisfies \( \err \lesssim \dimQ^{2}/n \).
Now the upper bound follows from Theorem~\ref{Tdevboundsharp}.
Similar arguments can be used for checking the lower bound by Theorem~\ref{Texpquadroi}.

%% file: qf-Gauss.tex

\Section{Moments of a Gaussian quadratic form}
\label{SmomentqfG}
Let \( \gaussv \) be standard normal in \( \R^{\dimp} \) for \( \dimp \leq \infty \).
Given a self-adjoint trace operator \( \BBH \), consider a quadratic form 
\( \bigl\langle \BBH \gaussv, \gaussv \bigr\rangle \).

\begin{lemma}
\label{Gaussmoments}
It holds \( \E \bigl\langle \BBH \gaussv, \gaussv \bigr\rangle = \tr \BBH \).
Moreover, 
\begin{EQA}
	\E \bigl( \bigl\langle \BBH \gaussv, \gaussv \bigr\rangle - \tr \BBH \bigr)^{2}
	&=&
	2 \tr \BBH^{2}  ,
	\\
	\E \bigl( \bigl\langle \BBH \gaussv, \gaussv \bigr\rangle - \tr \BBH \bigr)^{3}
	&=&
	8 \tr \BBH^{3} ,
	\\
	\E \bigl( \bigl\langle \BBH \gaussv, \gaussv \bigr\rangle - \tr \BBH \bigr)^{4}
	&=&
	48 \tr \BBH^{4} + 12 (\tr \BBH^{2})^{2} ,
	\\
	\E \bigl( \bigl\langle \BBH \gaussv, \gaussv \bigr\rangle - \tr \BBH \bigr)^{5}
	&=&
	512 \tr \BBH^{5} + 32 \tr \BBH^{2} \, \tr \BBH^{3} ,
\label{2pG2trD2DGm22m2}
\end{EQA}
and
\begin{EQA}
	\E \bigl\langle \BBH \gaussv, \gaussv \bigr\rangle^{2}
	&=&
	(\tr \BBH)^{2} + 2 \tr \BBH^{2},
	\\
	\E \bigl\langle \BBH \gaussv, \gaussv \bigr\rangle^{3}
	& = &
	(\tr \BBH)^{3} + 6 \tr \BBH \,\, \tr \BBH^{2} + 8 \tr \BBH^{3} ,
	\\
	\E \bigl\langle \BBH \gaussv, \gaussv \bigr\rangle^{4}
	& = &
	(\tr \BBH)^{4} + 12 (\tr \BBH)^{2} \tr \BBH^{2}
	+ 32 (\tr \BBH) \tr \BBH^{3}
	+ 48 \tr \BBH^{4} + 12 (\tr \BBH^{2})^{2} ,
\label{2pG2trD2DGm22m2}
	\\
	\Var \bigl\langle \BBH \gaussv, \gaussv \bigr\rangle^{2}
	& = &
	8 (\tr \BBH)^{2} \tr \BBH^{2}
	+ 32 (\tr \BBH) \tr \BBH^{3}
	+ 48 \tr \BBH^{4} + 8 (\tr \BBH^{2})^{2} .
\label{2pG2trD2DGm22m4}
\end{EQA}
Moreover, if \( \BBH \leq \Id_{\dimp} \) and \( \dimH = \tr \BBH \), then \( \tr \BBH^{m} \leq \dimH \| \BBH \|^{m-1} \) for 
\( m \geq 1 \) and
\begin{EQA}[rcccl]
	\E \bigl\langle \BBH \gaussv, \gaussv \bigr\rangle^{2}
	& \leq &
	\dimH^{2} + 2 \dimH \| \BBH \|
	&\leq &
	(\dimH + \| \BBH \|)^{2},
	\\
	\E \bigl\langle \BBH \gaussv, \gaussv \bigr\rangle^{3}
	& \leq &
	\dimH^{3} + 6 \dimH^{2} \| \BBH \| + 8 \dimH \| \BBH \|^{2}
	&\leq &
	(\dimH + 2 \| \BBH \|)^{3},
	\\
	\E \bigl\langle \BBH \gaussv, \gaussv \bigr\rangle^{4}
	& \leq &
	\dimH^{4} + 12 \dimH^{3} \| \BBH \|
	+ 44 \dimH^{2} \| \BBH \|^{2}
	+ 48 \dimH \| \BBH \|^{3}
	&\leq &
	(\dimH + 3 \| \BBH \|)^{4},
	\\
	\E \bigl\langle \BBH \gaussv, \gaussv \bigr\rangle^{5}
	& \leq &
	\dimH^{5} + 20 \dimH^{4} \| \BBH \| + 140 \dimH^{3} \| \BBH \|^{2} 
	+ 272 \dimH^{2} \| \BBH \|^{3} + 512 \dimH \| \BBH \|^{4}
	& \leq &
	(\dimH + 4 \| \BBH \|)^{5} \, .
\label{2pG2trD2DGm22m2}
	\\
	\Var \bigl\langle \BBH \gaussv, \gaussv \bigr\rangle^{2}
	& \leq &
	8 \dimH^{3} + 40 \dimH^{2} \| \BBH \| + 48 \dimH \| \BBH \|^{2}.
\label{2pG2trD2DGm22m4}
\end{EQA}
Finally,
\begin{EQA}
	\E (\gaussv \gaussv^{\T} - \Id_{\dimp}) \BBH (\gaussv \gaussv^{\T} - \Id_{\dimp}) 
	&=&
	\BBH + \tr (\BBH) \Id_{\dimp}
\label{njt66777888723fdgy}
\end{EQA}
yielding
\begin{EQA}
	\E \| \BBH (\gaussv \gaussv^{\T} - \Id_{\dimp}) \|_{\Fr}^{2}
	&=&
	(\tr \BBH)^{2} + \tr \BBH^{2} .
\label{njt66777888723fdgyf}
\end{EQA}
\end{lemma}

\begin{proof}
Let \( \gaussv \) be standard normal in \( \R^{\dimp} \).
The same holds for \( \Uv \gaussv \) for any orthogonal transform \( \Uv \) in \( \R^{\dimp} \).
The use of the spectral decomposition \( \BBH = \Uv^{\T} \Lambda \Uv \) with \( \Uv \) orthonormal and 
\( \Lambda \) diagonal enables us to represent 
\( \bigl\langle \BBH \gaussv, \gaussv \bigr\rangle = \bigl\langle \Lambda \Uv \gaussv, \Uv \gaussv \bigr\rangle \) and thus,
to reduce the statements to the case when \( \BBH \) is diagonal: \( \BBH = \diag(\supH_{1},\supH_{2},\ldots,\supH_{\dimp}) \).
Then 
\begin{EQA}
	\xi
	\eqdef
	\bigl\langle \BBH \gaussv, \gaussv \bigr\rangle - \tr \BBH
	&=&
	\sum_{j=1}^{\dimp} \supH_{j} (\gauss_{j}^{2} - 1) ,
\label{j1ljgj2m1}
\end{EQA}
where \( \gauss_{j} \) are i.i.d. standard normal. 
This easily yields with \( \dimH_{m} = \tr (\BBH^{m}) \)
\begin{EQA}
	\E \xi^{2}
	&=&
	\sum_{j=1}^{\dimp} \supH_{j}^{2} \E (\gauss_{j}^{2} - 1)^{2}
	=
	\E \chi^{2} \, \tr \BBH^{2} 
	=
	2 \dimH_{2} \, ,
	\\
	\E \xi^{3}
	&=&
	\sum_{j=1}^{\dimp} \supH_{j}^{3} \E (\gauss_{j}^{2} - 1)^{3}
	=
	\E \chi^{3} \, \tr \BBH^{3} 
	=
	8 \dimH_{3} \, ,
	\\
	\E \xi^{4}
	&=&
	\sum_{j=1}^{\dimp} \supH_{j}^{4} (\gauss_{j}^{2} - 1)^{4}
	+ \sum_{i\neq j} \supH_{i}^{2} \supH_{j}^{2} \E (\gauss_{i}^{2} - 1)^{2} \E (\gauss_{j}^{2} - 1)^{2}
	\\
	&=&
	\bigl( \E \chi^{4} - 3 (\E \chi^{2})^{2} \bigr) \tr \BBH^{4} + 3 (\E \chi^{2} \, \tr \BBH^{2})^{2}
	=
	48 \dimH_{4} + 12 \dimH_{2}^{2} ,
	\\
	\E \xi^{5}
	&=&
	\sum_{j=1}^{\dimp} \supH_{j}^{5} (\gauss_{j}^{2} - 1)^{5}
	+ \sum_{i\neq j} \supH_{i}^{2} \supH_{j}^{3} \E (\gauss_{i}^{2} - 1)^{2} \E (\gauss_{j}^{2} - 1)^{3}
	\\
	&=&
	\bigl\{ \E (\gauss^{2} - 1)^{5} - \E (\gauss^{2} - 1)^{2} \E (\gauss^{2} - 1)^{3} \bigr\} \, \tr \BBH^{5} 
	+ \E (\gauss^{2} - 1)^{2} \E (\gauss^{2} - 1)^{3} \, \tr \BBH^{2} \, \tr \BBH^{3}
	\\
	&=&
	512 \dimH_{5} + 32 \dimH_{2} \, \dimH_{3} \, .
\label{2pG2trD2DGm22m2}
\end{EQA}
and
\begin{EQA}
	\E \bigl\langle \BBH \gaussv, \gaussv \bigr\rangle^{2}
	&=&
	\bigl( \E \bigl\langle \BBH \gaussv, \gaussv \bigr\rangle \bigr)^{2} 
	+ \E \xi^{2}
	= 
	\dimH^{2} + 2 \dimH_{2} \, ,
	\\
	\E \bigl\langle \BBH \gaussv, \gaussv \bigr\rangle^{3}
	& = &
	\E ( \xi + \dimH )^{3}
	=
	\dimH^{3} + \E \xi^{3} + 3 \dimH \,\, \E \xi^{2}
	=
	\dimH^{3} + 6 \dimH \,\, \dimH_{2} + 8 \dimH_{3} ,
	\\
	\E \bigl\langle \BBH \gaussv, \gaussv \bigr\rangle^{4}
	& = &
	\E \bigl( \xi + \dimH \bigr)^{4}
	=
	\dimH^{4} + 6 \dimH^{2} \E \xi^{2} + 4 \dimH \, \E \xi^{3} + \E \xi^{4}
	\\
	&=&
	\dimH^{4} + 12 \dimH^{2} \, \dimH_{2}
	+ 32 \dimH \, \dimH_{3}
	+ 48 \dimH_{4} + 12 \dimH_{2}^{2} ,
\label{2pG2trD2DGm22m2}
\end{EQA}
and 
\begin{EQA}
	&& \nquad
	\Var \bigl\langle \BBH \gaussv, \gaussv \bigr\rangle^{2}
	= 
	\E ( \xi + \dimH )^{4}
	- \bigl( \dimH^{2} + 2 \dimH_{2} \bigr)^{2}
	\\
	&=&
	\dimH^{4} + 6 \dimH^{2} \E \xi^{2} + 4 \dimH \, \E \xi^{3} + \E \xi^{4}
	- \bigl( \dimH^{2} + 2 \dimH_{2} \bigr)^{2}
	= 
	8 \dimH^{2} \, \dimH_{2}
	+ 32 \dimH \, \dimH_{3}
	+ 48 \dimH_{4} + 8 \dimH_{2}^{2} \, .
\label{2pG2trD2DGm22m4}
\end{EQA}
Also
\begin{EQA}
	\E \bigl\langle \BBH \gaussv, \gaussv \bigr\rangle^{5}
	& = &
	\E \bigl( \xi + \dimH \bigr)^{5}
	=
	\dimH^{5} + 10 \dimH^{3} \E \xi^{2} + 10 \dimH^{2} \, \E \xi^{3} + 5 \dimH \E \xi^{4}
	+ \E \xi^{5}
	\\
	&=&
	\dimH^{5} + 20 \dimH^{3} \, \dimH_{2} + 80 \dimH^{2} \dimH_{3} 
	+ 5 \dimH (48 \dimH_{4} + 12 \dimH_{2}^{2})
	+ 512 \dimH_{5} + 32 \dimH_{2} \, \dimH_{3} \, .
\label{2pG2trD2DGm22m25}
\end{EQA}
Assume \( \| \BBH \| = 1 \) yielding \( \dimH_{m} \leq \dimH \).
Then
\begin{EQA}
	\E \bigl\langle \BBH \gaussv, \gaussv \bigr\rangle^{2}
	& \leq &
	\dimH^{2} + 2 \dimH
	\leq 
	(\dimH + 1)^{2} \, ,
	\\
	\E \bigl\langle \BBH \gaussv, \gaussv \bigr\rangle^{3} 
	& \leq &
	\dimH^{3} + 6 \dimH^{2} + 8 \dimH
	\leq 
	(\dimH + 2)^{3} ,
	\\
	\E \bigl\langle \BBH \gaussv, \gaussv \bigr\rangle^{4}
	& \leq &
	\dimH^{4} + 12 \dimH^{3} + 44 \dimH^{2} + 48 \dimH 
	\leq 
	(\dimH + 3)^{4} ,
	\\
	\E \bigl\langle \BBH \gaussv, \gaussv \bigr\rangle^{5}
	& \leq &
	\dimH^{5} + 20 \dimH^{4} + 140 \dimH^{3} + 272 \dimH^{2} + 512 \dimH 
	\leq 
	(\dimH + 4)^{5} \, .
\label{2pG2trD2DGm22m25}
\end{EQA}
For the last result of the lemma, observe that with \( \BBH = \diag(\supH_{1},\supH_{2},\ldots,\supH_{\dimp}) \), 
\begin{EQA}
	\E \| \BBH^{1/2} (\gaussv \gaussv^{\T} - \Id_{\dimp}) \BBH^{1/2} \|_{\Fr}^{2}
	&=&
	\sum_{i,j=1}^{\dimp} \supH_{i} \supH_{j} \E (\gauss_{i} \gauss_{j} - \delta_{i,j})^{2}
	=
	\left( \sum_{i=1}^{\dimp} \supH_{i} \right)^{2} + \sum_{i=1}^{\dimp} \supH_{i}^{2} 
\label{ikcduywjwsiv98emdvuw}
\end{EQA}
and assertion \eqref{njt66777888723fdgyf} follows.
\end{proof}

Now we compute the exponential moments of centered and non-centered quadratic forms.

\begin{lemma}
\label{Lqfexpmom}
Let \( \| \BBH \| = \supH \) and \( \gaussv \sim \ND(0,\Id_{\dimp}) \).
Then for any \( \mu \in (0,\supH^{-1}) \), 
\begin{EQA}
	\E \exp \Bigl\{ \frac{\mu}{2} \langle \BBH \gaussv, \gaussv \rangle \Bigr\}
	&=&
	\det(\Id_{\dimp} - \mu \BBH)^{-1/2} \, .
\label{m2v241m41m}
\end{EQA}
Moreover, with \( \dimH = \tr \BBH \) and \( \vH^{2} = \tr \BBH^{2} \)
\begin{EQA}
	\log \E \exp \Bigl\{ \frac{\mu}{2} \bigl( \langle \BBH \gaussv, \gaussv \rangle - \dimH \bigr) \Bigr\}
	& \leq &
	\frac{\mu^{2} \vH^{2}}{4 (1 - \supH \mu)} \, .
\label{m2v241m41mb}
\end{EQA}
If \( \BBH \) is positive semidefinite, \( \supH_{j} \geq 0 \), then 
\begin{EQA}
	\log \E \exp \Bigl\{ - \frac{\mu}{2} \bigl( \langle \BBH \gaussv, \gaussv \rangle - \dimH \bigr) \Bigr\}
	& \leq &
	\frac{\mu^{2} \vH^{2}}{4} \, .
\label{m2v241m41mbn}
\end{EQA}
For any complex valued \( \muH \) with \( \supH |\muH| < 1 \),
\begin{EQA}
	\biggl| \log \E \exp \Bigl\{ 
			\frac{\mu}{2} \bigl( \langle \BBH \gaussv, \gaussv \rangle - \dimH \bigr) - \frac{\muH^{2} \tr \BBH^{2}}{4}
		\Bigr\} 
	\biggr|
	& \leq &
	\frac{\supH |\mu|^{3} \vH^{2} }{6 (1 - \supH |\mu|)} \, .
\label{vbu7j3hg8hryhghyidegwdg}
\end{EQA}
\end{lemma}

\begin{proof}
W.l.o.g. assume \( \supH = 1 \).
Let \( \supH_{j} \) be the eigenvalues of \( \BBH \), \( |\supH_{j}| \leq 1 \).
As in Lemma~\ref{Gaussmoments}, one can reduce the statement to the case of a diagonal matrix 
\( \BBH = \diag\bigl( \supH_{j} \bigr) \). 
Then \( \langle \BBH \gaussv, \gaussv \rangle = \sum_{j=1}^{\dimp} \supH_{j} \gauss_{j}^{2} \) and 
by independence of the \( \gauss_{j} \)'s
\begin{EQA}
	&& \nquad
	\E \Bigl\{ \frac{\mu}{2} \langle \BBH \gaussv, \gaussv \rangle  \Bigr\}
	=
	\prod_{j=1}^{\dimp} \E \exp \Bigl( \frac{\mu}{2} \supH_{j} \eps_{j}^{2} \Bigr)
	=
	\prod_{j=1}^{\dimp} \frac{1}{\sqrt{1 - \mu \supH_{j}}} 
	=
	\det \bigl( \Id_{\dimp} - \mu \BBH \bigr)^{-1/2} .
\label{dOImuBm12EB}
\end{EQA}
Below we use the simple bounds: 
\begin{EQ}[rcl]
\label{lo1uusk2iukkp}
	- \log(1 - u) - u
	&=&
	\sum_{k=2}^{\infty} \frac{u^{k}}{k}
	\leq 
	\frac{u^{2}}{2} \sum_{k=0}^{\infty} u^{k} 
	=
	\frac{u^{2}}{2 (1 - u)} \, ,
	\qquad 
	u \in (0,1),
	\\
	- \log(1 - u) + u
	&=&
	\sum_{k=2}^{\infty} \frac{u^{k}}{k}
	\leq 
	\frac{u^{2}}{2} \, ,
	\qquad \qquad
	u \in (-1,0).
\label{lo1uusk2iukk}
\end{EQ}
Now it holds for \( \mu > 0 \)
\begin{EQA}
	&& \nquad
	\log \E \Bigl\{ \frac{\mu}{2} \bigl( \langle \BBH \gaussv, \gaussv \rangle - \dimH \bigr) \Bigr\}
	=
	\log \det(\Id_{\dimp} - \mu \BBH)^{-1/2} - \frac{\mu \, \dimH}{2}
	\\
	&=&
	- \frac{1}{2} \sum_{j=1}^{\dimp} \bigl\{ \log(1 - \mu \supH_{j}) + \mu \supH_{j} \bigr\}
	\leq 
	\sum_{j=1}^{\dimp} \frac{\mu^{2} \supH_{j}^{2}}{4 (1 - \mu \supH_{j})} 
	\leq 
	\frac{\mu^{2} \vH^{2}}{4 (1 - \mu \supH)} \, .
\label{m2v241m4mj1pd}
\end{EQA}
Similarly for any complex \( \mu \) with \( |\mu| \supH < 1 \)
\begin{EQA}
	&& \nquad
	\left| 
		\log \E \Bigl\{ \frac{\mu}{2} \bigl( \langle \BBH \gaussv, \gaussv \rangle - \dimH \bigr) 
		- \frac{\muH^{2} \tr \BBH^{2}}{4} \Bigr\}
	\right|
	=
	\left| \log \det(\Id_{\dimp} - \mu \BBH)^{-1/2} - \frac{\mu \, \dimH}{2} - \frac{\muH^{2} \tr \BBH^{2}}{4} \right|
	\\
	&=&
	\frac{1}{2} \left| 
		\sum_{j=1}^{\dimp} \biggl\{ \log(1 - \mu \supH_{j}) - \mu \supH_{j} - \frac{\mu^{2} \supH_{j}^{2}}{2} \biggr\} 
	\right|
	\leq 
	\sum_{j=1}^{\dimp} \frac{|\mu \supH_{j}|^{3}}{6 (1 - \supH |\mu|)} 
	=
	\frac{| \mu |^{3} \supH \vH^{2}}{6 (1 - \supH |\mu|)} \, .
\label{m2v241m4mj1pd}
\end{EQA}
Statement \eqref{m2v241m41mbn} can be proved similarly.
\end{proof}

Now we consider the case of a non-centered quadratic form
\( \langle \BBH \gaussv,\gaussv \rangle/2 + \langle \Av,\gaussv \rangle \) for a fixed vector \( \Av \).

\begin{lemma}
\label{Lexpmomnoncen}
Let \( \| \BBH \| = \supH < 1 \). 
Then for any \( \Av \)
\begin{EQA}
	\E \exp\Bigl\{ \frac{1}{2}\langle \BBH \gaussv,\gaussv \rangle + \langle \Av,\gaussv \rangle \Bigr\}
	&=&
	\exp\Bigl\{ \frac{\| (\Id_{\dimp} - \BBH)^{-1/2} \Av \|^{2}}{2} \Bigr\} \, \det(\Id_{\dimp} - \BBH)^{-1/2} .
\label{EeBf12BggA}
\end{EQA}
Moreover, for any \( \mu \in (0,1) \)
\begin{EQA}
	&& \nquad
	\log \E \exp\Bigl\{ 
		\frac{\mu}{2} \bigl( \langle \BBH \gaussv,\gaussv \rangle - \dimH \bigr) + \langle \Av,\gaussv \rangle 
	\Bigr\}
	\\
	&=&
	\frac{\| (\Id_{\dimp} - \mu \BBH)^{-1/2} \Av \|^{2}}{2} + \log \det(\Id_{\dimp} - \mu \BBH)^{-1/2} - \mu \, \dimH 
	\\
	& \leq &
	\frac{\| (\Id_{\dimp} - \mu \BBH)^{-1/2} \Av \|^{2}}{2} + \frac{\mu^{2} \vH^{2}}{4 (1 - \supH \mu)} \, .
\label{EeBf12BggAmu}
\end{EQA}
\end{lemma}

\begin{proof}
Denote \( \av = (\Id_{\dimp} - \BBH)^{-1/2} \Av \). 
It holds by change of variables \( (\Id_{\dimp} - \BBH)^{1/2} \xv = \uv \) for \( \CONSTi_{\dimp} = (2\pi)^{-\dimp/2} \)
\begin{EQA}
	&& \nquad
	\E \exp\Bigl\{ \frac{1}{2}\langle \BBH \gaussv,\gaussv \rangle + \langle \Av,\gaussv \rangle \Bigr\}
	=
	\CONSTi_{\dimp}
	\int \exp\Bigl\{ - \frac{1}{2}\langle (\Id_{\dimp} - \BBH) \xv,\xv \rangle + \langle \Av,\xv \rangle \Bigr\} d\xv
	\\
	&=&
	\CONSTi_{\dimp}
	\det(\Id_{\dimp} - \BBH)^{-1/2}
	\int \exp\Bigl\{ - \frac{1}{2} \| \uv \|^{2} + \langle \av,\uv \rangle \Bigr\} d\uv
	=
	\det(\Id_{\dimp} - \BBH)^{-1/2} \, 	\ex^{\| \av \|^{2}/2}  	.
\label{EeBf12BggAp}
\end{EQA}
The last inequality \eqref{EeBf12BggAmu} follows by \eqref{m2v241m41mb}.
\end{proof}

\Section{Deviation bounds for Gaussian quadratic forms}
\label{SdevboundGauss}
The next result explains the concentration effect of \( \| \QP \xiv \|^{2} \)
for a centered Gaussian vector \( \xiv \sim \ND(0,\HVB^{2}) \) and a linear operator \( \QP \colon \R^{\dimp} \to \R^{\dimq} \),
\( \dimp,\dimq \leq \infty \).
We use a version from \cite{laurentmassart2000}.
For completeness, we present a simple proof.

\begin{theorem}
\label{TexpbLGA}
\label{Lxiv2LD}
\label{Cuvepsuv0}
Let \( \xiv \sim \ND(0,\HVB^{2}) \) be a Gaussian element in \( \R^{\dimp} \) and let
\( \QP \colon \R^{\dimp} \to \R^{\dimq} \) be such that \( \BBH = \QP \HVB^{2} \QP^{\T} \) 
is a trace operator in \( \R^{\dimq} \).
Then with \( \dimH = \tr(\BBH) \), \( \vH^{2} = \tr(\BBH^{2}) \), and \( \supH = \| \BBH \| \),
it holds for any \( \xx \geq 0 \)
\begin{EQA}
\label{Pxiv2dimAvp12}
	\P\Bigl( \| \QP \xiv \|^{2} - \dimH > 2 \vH \, \sqrt{\xx} + 2 \supH \xx \Bigr)
	& \leq &
	\ex^{-\xx} ,
	\\
	\P\Bigl( \| \QP \xiv \|^{2} - \dimH \leq - 2 \vH \, \sqrt{\xx} \Bigr)
	& \leq &
	\ex^{-\xx} .
\label{Pxiv2dimAvp12m}
\end{EQA}
It also implies 
\begin{EQA}
	\P\bigl( \bigl| \| \QP \xiv \|^{2} - \dimH \bigr| > \zq_{2}(\BBH,\xx) \bigr)
	& \leq &
	2 \ex^{-\xx} ,
\label{PxivTBBdimA2vp}
\end{EQA}
with
\begin{EQA}
	\zq_{2}(\BBH,\xx)
	& \eqdef &
	2 \vH \, \sqrt{\xx} + 2 \supH \xx \,\, .
\label{zqdefGQF}
\end{EQA}
%
\end{theorem}

\begin{proof}
We use the identity in distribution \( \| \QP \xiv \|^{2} \eqd \langle \BBH \gaussv, \gaussv \rangle \) with
 \( \gaussv \sim \ND(0,\Id_{\dimq}) \).
Markov's inequality yields for any \( \mu > 0 \)
\begin{EQA}
	\P\Bigl( \langle \BBH \gaussv, \gaussv \rangle - \dimH > \zq_{2}(\BBH,\xx) \Bigr)
	& \leq &
	\E \exp \Bigl( \frac{\mu}{2} \bigl( \langle \BBH \gaussv, \gaussv \rangle - \dimH \bigr) - \frac{\mu \, \zq_{2}(\BBH,\xx)}{2} 
	\Bigr) \, .
\label{PBggiz2E2mz2}
\end{EQA}
Given \( \xx > 0 \), fix \( \mu < 1/\supH \) by the equation
\begin{EQA}
	\frac{\mu}{1 - \supH \mu} 
	&=&
	\frac{2 \sqrt{\xx}}{\vH} \, 
	\quad \text{ or } \quad
	\mu^{-1} 
	=
	\supH + \frac{\vH}{2 \sqrt{\xx}} \, .
\label{1v2sxm12m1m}
\end{EQA}
By \eqref{m2v241m41mb}
\begin{EQA}
	&& \nquad
	\log \E \Bigl\{ \frac{\mu}{2} \bigl( \langle \BBH \gaussv, \gaussv \rangle - \dimH \bigr) \Bigr\}
	\leq 
	\frac{\mu^{2} \vH^{2}}{4 (1 - \supH \mu)} \, .
\label{m2v241m4mj1p}
\end{EQA}
For \eqref{Pxiv2dimAvp12}, it remains to check that the choice \( \mu \) by \eqref{1v2sxm12m1m} yields
\begin{EQA}
	\frac{\mu^{2} \vH^{2}}{4 (1 - \supH \mu)} - \frac{\mu \, \zq_{2}(\BBH,\xx)}{2}
	& = &
	\frac{\mu^{2} \vH^{2}}{4 (1 - \supH \mu)} - \mu \bigl( \vH \sqrt{\xx} + \supH \xx \bigr)
	=
	\mu \Bigl( \frac{\vH \sqrt{\xx}}{2} - \vH \sqrt{\xx} - \supH \xx \Bigr)
	=
	- \xx .
\label{m2vA241muz2}
\end{EQA}
The bound \eqref{Pxiv2dimAvp12m} is obtained similarly from Markov's inequality 
applied to \( - \langle \BBH \gaussv, \gaussv \rangle + \dimH \) with \( \mu = 2 \vH^{-1} \sqrt{\xx} \).
The use of \eqref{m2v241m41mbn} yields
\begin{EQA}
	&& \nquad
	\P\Bigl( \langle \BBH \gaussv, \gaussv \rangle - \dimH < - 2 \vH \sqrt{\xx} \Bigr)
	\leq
	\E \exp \Bigl\{ \frac{\mu}{2} \bigl( - \langle \BBH \gaussv, \gaussv \rangle + \dimH \bigr) - \mu \, \vH \sqrt{\xx} 
	\Bigr\}
	\\
	& \leq &
	\exp \Bigl( \frac{\mu^{2} \vH^{2}}{4} - \mu \, \vH \sqrt{\xx} \Bigr) 
	=
	\ex^{-\xx} \, 
\label{PBggiz2E2mz2}
\end{EQA}
as required.
\end{proof}

\begin{corollary}
\label{CTexpbLGAd}
Assume the conditions of Theorem~\ref{TexpbLGA}.
Then for \( \zq > \vH \)
\begin{EQA}
	\P\bigl( \bigl| \| \QP \xiv \|^{2} - \dimH \bigr| \ge \zq \bigr)
	& \leq &
	2 \exp\biggl\{ - \frac{\zq^{2}}{\bigl( \vH + \sqrt{\vH^{2} + 2 \supH \zq} \bigr)^{2}} \biggr\}
	\leq 
	2 \exp\biggl( - \frac{\zq^{2}}{4\vH^{2} + 4 \supH \zq} \biggr) .
	\qquad
	\qquad
\label{3z2spsp2z3z2}
\end{EQA}
\end{corollary}

\begin{proof}
Given \( \zq \), define \( \xx \) by 
\( 2 \vH \sqrt{\xx} + 2 \supH \xx = \zq \) or 
\( 2 \supH \sqrt{\xx} = \sqrt{\vH^{2} + 2 \supH \zq} - \vH \).
Then
\begin{EQA}
	\P\bigl( \| \QP \xiv \|^{2} - \dimH \ge \zq \bigr)
	& \leq &
	\ex^{-\xx} 
	=
	\exp\biggl\{ - \frac{\bigl( \sqrt{\vH^{2} + 2 \supH \zq} - \vH \bigr)^{2}}{4 \supH^{2}} \biggr\}
	=
	\exp\biggl\{ - \frac{\zq^{2}}{\bigl( \vH + \sqrt{\vH^{2} + 2 \supH \zq} \bigr)^{2}} \biggr\}.
\label{3emzmsp22z2c}
\end{EQA}
This yields \eqref{3z2spsp2z3z2} by direct calculus.
\end{proof}

Of course, bound \eqref{3z2spsp2z3z2} is sensible only if \( \zq \gg \vH \).

\ifadap{}{
\begin{corollary}
\label{RsochpHsA}
Assume the conditions of Theorem~\ref{TexpbLGA}.
If also \( \BBH \geq 0 \), then 
\begin{EQA}
\label{Pxiv2dimAxx12}
	\P\Bigl( \| \QP \xiv \|^{2} \geq \zq^{2}(\BBH,\xx) \Bigr)
	& \leq &
	\ex^{-\xx} 
\end{EQA}
with 
\begin{EQA}
	\zq^{2}(\BBH,\xx)
	& \eqdef &
	\dimH + 2 \vH \, \sqrt{\xx} + 2 \supH \xx
	\leq 
	\bigl( \sqrt{\dimH} + \sqrt{2 \supH \xx} \bigr)^{2} \, .
\label{zzxxppdBlroBB}
\end{EQA}
Also
\begin{EQA}
	\P\Bigl( \| \QP \xiv \|^{2} - \dimH < - 2 \vH \, \sqrt{\xx} \Bigr)
	& \leq &
	\ex^{-\xx} .
\label{Pxiv2dimAvp12d}
\end{EQA}
\end{corollary}

\begin{proof}
The definition implies \( \vH^{2} \leq \dimH \supH \)
yielding the statement of the corollary.
\end{proof}

As a special case, we present a bound for the chi-squared distribution 
corresponding to \( \QP = \HVB^{2} = \Id_{\dimp} \), \( \dimp < \infty \).
Then \( \BBH = \Id_{\dimp} \), \( \tr (\BBH) = \dimp \), \( \tr(\BBH^{2}) = \dimp \) and \( \supH(\BBH) = 1 \).

\begin{corollary}
\label{Cchi2p}
Let \( \gaussv \) be a standard normal vector in \( \R^{\dimp} \).
Then for any \( \xx > 0 \)
\begin{EQA}[ccl]
\label{Pxi2pm2px}
	\P\bigl( \| \gaussv \|^{2} \geq \dimp + 2 \sqrt{\dimp \, \xx} + 2 \xx \bigr)
	& \leq &
	\ex^{-\xx},
	\\
	\P\bigl( \| \gaussv \| \,\,  \geq \sqrt{\dimp} + \sqrt{2 \xx} \bigr)
	& \leq &
	\ex^{-\xx} ,
\label{Pxi2pm2px12}
	\\
	\P\bigl( \| \gaussv \|^{2} \leq \dimp - 2 \sqrt{\dimp \, \xx} \bigr)
	& \leq &
	\ex^{-\xx}	.
\label{Pxi2pm2px22}
\end{EQA}
\end{corollary}
}

The bound of Theorem~\ref{TexpbLGA} 
can be represented as a usual deviation bound.

\begin{theorem}
\label{CTexpbLGA}
Assume the conditions of Theorem~\ref{TexpbLGA}.
For \( \yy > 0 \), define
\begin{EQA}
	\xx(\yy)
	& \eqdef &
	\frac{(\sqrt{\yy + \dimH} - \sqrt{\dimH})^{2}}{4 \supH} \, .
\label{iuvfiiow3kboieheuf}
\end{EQA}
Then
\begin{EQA}
	\P\bigl( \| \QP \xiv \|^{2} \ge \dimH + \yy \bigr)
	& \leq &
	\ex^{- \xx(\yy)} ,
\label{3emzmsp22z2}
	\\
	\E \bigl\{ (\| \QP \xiv \|^{2} - \dimH) \Ind\bigl( \| \QP \xiv \|^{2} \ge \dimH + \yy \bigr) \bigr\}
	& \leq &
	2 \Bigl( \frac{\yy + \dimH}{\supH \, \xx(\yy)} \Bigr)^{1/2} \, \, 
	\ex^{- \xx(\yy)} \, .
	\qquad
	\quad
\label{3emzmsp22z2e}
\end{EQA}
Moreover, let \( \muH > 0 \) fulfill \( \rexH =  \muH \supH + \muH \sqrt{\supH \dimH / \xx(\yy)} < 1 \). 
Then 
\begin{EQA}
	\E \bigl\{ \ex^{\muH (\| \QP \xiv \|^{2} - \dimH)/2} \Ind( \| \QP \xiv \|^{2} \ge \dimH + \yy) \bigr\}
	& \leq &
	\frac{1}{1 - \rexH} \, \exp\{ - (1 - \rexH) \xx(\yy) \} \, .
	\qquad
\label{llkknbononjm9hig4e}
\end{EQA}
\end{theorem}

\begin{proof}
Normalizing by \( \supH \) reduces the statements to the case with \( \supH = 1 \).
Define \( \eta = \| \QP \xiv \|^{2} - \dimH \) 
and
\begin{EQA}
	\zq(\xx)
	&=&
	2 \sqrt{\dimH \, \xx} + 2 \xx .
\label{0kmuy765433udgswhhh}
\end{EQA}
Then by \eqref{Pxiv2dimAvp12} \( \P(\eta \geq \zq(\xx)) \leq \ex^{-\xx} \).
Inverting the relation \eqref{0kmuy765433udgswhhh} yields
\begin{EQA}
	\xx(\zq)
	&=&
	\frac{1}{4} \bigl( \sqrt{\zq + \dimH} - \sqrt{\dimH} \bigr)^{2}
\label{jkv78fdjryfgsdfghgj}
\end{EQA}
and \eqref{3emzmsp22z2} follows by applying \( \zq = \yy \).
Further, 
\begin{EQA}
	\E \bigl\{ \eta \Ind(\eta \geq \yy) \bigr\}
	&=&
	\int_{\yy}^{\infty} \P(\eta \geq \zq) \, d\zq
	\leq 
	\int_{\yy}^{\infty} \ex^{ - \xx(\zq) } \, d\zq
	= 
	\int_{\xx(\yy)}^{\infty} \ex^{-\xx} \, \zq'(\xx) \, d\xx \, .
\label{zEe2Iezz2c2H23}
\end{EQA} 
As \( \zq'(\xx) = 2 + \sqrt{\dimH/\xx} \) monotonously decreases with \( \xx \), we derive
\begin{EQA}
	\E \bigl\{ \eta \Ind(\eta \geq \yy) \bigr\}
	& \leq &
	\zq'(\xx(\yy)) \ex^{-\xx(\yy)}
	=
	\frac{1}{\xx'(\yy)} \, \ex^{- \xx(\yy)}
	=
	\frac{4 \sqrt{\yy + \dimH}}{\sqrt{\yy + \dimH} - \sqrt{\dimH}} \, \ex^{- \xx(\yy)}
\label{e7ygv76bgughytuj}
\end{EQA}
and \eqref{3emzmsp22z2e} follows.

In a similar way, define \( \zqe(\xx) \) from the relation
\( \muH^{-1} \log \zqe(\xx) = \sqrt{\dimH \, \xx} + \xx \) yielding
\begin{EQA}
	\zqe(\xx)
	&=&
	\exp \bigl( \muH \sqrt{\dimH \, \xx} + \muH \, \xx \bigr) .
\label{jvcjjuvue37r6gtur4r}
\end{EQA}
The inverse relation reads
\begin{EQA}
	\xxe(\zqe)
	&=&
	\bigl( \sqrt{\muH^{-1} \log \zqe + \dimH/4} - \sqrt{\dimH/4} \bigr)^{2} .
\label{jkv78fdjryfgsdfghgjex}
\end{EQA}
Then with \( \xx(\yy) = \xxe(\ex^{\muH \yy/2}) = \bigl( \sqrt{\yy + \dimH} - \sqrt{\dimH} \bigr)^{2}/4 \)
\begin{EQA}
	\E \bigl\{ \ex^{\muH \eta/2} \Ind(\eta \geq \yy) \bigr\}
	&=&
	\int_{\ex^{\muH \yy/2}}^{\infty} \P(\ex^{\muH \eta/2} \geq \zqe) \, d\zqe
	=
	\int_{\ex^{\muH \yy/2}}^{\infty} \P(\eta \geq 2\muH^{-1} \log \zqe) \, d\zqe
	\\
	& \leq &
	\int_{\ex^{\muH \yy/2}}^{\infty} \ex^{ - \xxe(\zqe) } \, d\zqe
	= 
	\int_{\xx(\yy)}^{\infty} \ex^{-\xx} \, \zqe'(\xx) \, d\xx .
\label{zEe2Iezz2c2H23}
\end{EQA} 
Further, in view of \( \muH + 0.5 \,\muH \sqrt{\dimH/\xx} < \muH + \muH \sqrt{\dimH / \xx(\yy)} = \rexH < 1 \) for 
\( \xx \geq \xx(\yy) \), it holds
\begin{EQA}
	\zqe'(\xx)
	&=&
	\bigl( \muH + 0.5 \, \muH \sqrt{\dimH/\xx} \bigr) \exp \bigl( \muH \sqrt{\dimH \, \xx} + \muH \, \xx \bigr) 
	\leq 
	\exp \bigl( \muH \, \xx \sqrt{\dimH / \xx(\yy)} + \muH \, \xx \bigr)
	=
	\exp (\rexH \, \xx) 
\label{jcuyu3ww3jbkihjitwedk}
\end{EQA}
and  
\begin{EQA}
	\E \bigl\{ \ex^{\muH \eta/2} \Ind(\eta \geq \yy) \bigr\}
	& \leq &
	\int_{\xx(\yy)}^{\infty} \ex^{-(1 - \rexH)\xx} \, d\xx 
	=
	\frac{1}{1 - \rexH} \, \ex^{-(1 - \rexH)\xx(\yy)} \, 
\label{zEe2Iezz2c2H23}
\end{EQA} 
and \eqref{llkknbononjm9hig4e} follows.
\end{proof}

%% file: t3-Gauss4.tex
\Section{Some bounds for a third order Gaussian tensor}
\label{SDB3tens}
Let \( \Tens = \bigl( \Tens_{i,j,k} \bigr) \) be a third order symmetric tensor, that is,
\( \Tens_{i,j,k} = \Tens_{\perm(i,j,k)} \) for any permutation \( \perm \) of the triple \( (i,j,k) \).
%
This section present a deviation bound for a Gaussian tensor sum 
\( \Tens(\GaussD) \eqdef \langle \Tens, \GaussD^{\otimes 3} \rangle \) 
for a Gaussian zero mean vector \( \GaussD \sim \ND(0,\DPTG^{-2}) \) in \( \R^{\dimp} \).
Much more general results for higher order tensors are available in the literature, see e.g. \cite{GSS2021} and \cite{AW2013} and references therein.
We, however, present an independent self-contained study which delivers finite sample and sharp results.
Later we use notations
\begin{EQA}[ccl]
	\| \Tens \|
	&=&
	\sup_{\| \uv_{1} \| = \| \uv_{2} \| = \| \uv_{3} \| = 1} 
		\bigl| \langle \Tens, \uv_{1} \otimes \uv_{2} \otimes \uv_{3} \rangle \bigr| \, .
\label{jcu7ejufd76ehyew236h}
\end{EQA}
Banach's characterization \cite{Banach1938,nie2017} yields 
\begin{EQA}
	\| \Tens \|
	&=&
	\sup_{\| \uv \| = 1} \bigl| \langle \Tens, \uv^{\otimes 3} \rangle \bigr| \, . 
\label{iv7dcu3921kfvhw7fjew}
\end{EQA}
Define
\begin{EQA}
	\Tens(\uv)
	&=&
	\langle \Tens, \uv^{\otimes 3} \rangle
	=
	\sum_{i,j,k=1}^{\dimp} \Tens_{i,j,k} \, u_{i} \, u_{j} \, u_{k} \, ,
	\qquad
	\uv = (u_{i}) \in \R^{\dimp} \, .
\label{vhjr8vh8r84423rghyrf}
\end{EQA}
Clearly \( \Tens(\uv) \) is a third order polynomial function on \( \R^{\dimp} \).
Define also its gradient \( \nabla \Tens(\uv) \in \R^{\dimp} \).
Each entry of the gradient vector \( \nabla \Tens(\uv) \) is a second order polynomial of \( \uv \).
Symmetricity of \( \Tens \) implies for any \( \uv \in \R^{\dimp} \)
\begin{EQ}[rcccl]
	\nabla \Tens(\uv)
	&=&
	\biggl( 3 \sum_{j,k=1}^{\dimp} \Tens_{i,j,k} \, u_{j} \, u_{k}
	\biggr)_{i=1,\ldots,\dimp} 
	&=&
	3 \bigl( \langle \Tens_{i} \, , \uv \otimes \uv \rangle \bigr)_{i=1,\ldots,\dimp} \,\, ,
	\\
	\nabla^{2} \Tens(\uv)
	&=&
	\biggl( 6 \sum_{i=1}^{\dimp} \Tens_{i,j,k} \, u_{i} 
	\biggr)_{j,k=1,\ldots,\dimp} 
	&=&
	6 \Tens[\uv] \, ,
\label{ho5r9hkrh4ygj32fjur}
\end{EQ}
where \( \Tens_{i} \) is the sub-tensor of order 2 with \( (\Tens_{i})_{j,k} = \Tens_{i,j,k} \) and
\begin{EQA}
	\Tens[\uv]
	& \eqdef &
	\sum_{i=1}^{\dimp} u_{i} \, \Tens_{i} \, .
\label{ufhjdeuvbuftjeee3jrjc3}
\end{EQA} 
Also
\begin{EQA}
	\Tens(\uv)
	&=&
	\frac{1}{3} \langle \nabla \Tens(\uv),\uv \rangle
	=
	\frac{1}{6} \langle \nabla^{2} \Tens(\uv),\uv^{\otimes 2} \rangle \, .
\label{tv8f8e35etfghr67ghj}
\end{EQA}
%
For the norm of the vector \( \nabla \Tens(\uv) \) and of the matrix \( \nabla^{2}\Tens(\uv) \), it holds by \eqref{jcu7ejufd76ehyew236h}
\begin{EQA}
	\| \nabla \Tens(\uv) \|
	&=&
	\sup_{\omegav \in \R^{\dimp} \colon \| \omegav \| = 1} \langle \nabla \Tens(\uv),\omegav \rangle
	=
	\sup_{\omegav \in \R^{\dimp} \colon \| \omegav \| = 1} 3 \langle \Tens,\uv \otimes \uv \otimes \omegav \rangle
	=
	3 \| \Tens \| \, \| \uv \|^{2},
	\\
	\| \nabla^{2} \Tens(\uv) \|
	&=&
	\sup_{\omegav \in \R^{\dimp} \colon \| \omegav \| = 1} 
	\bigl| \bigl\langle \nabla^{2} \Tens(\uv) \, , \omegav \otimes \omegav \bigr\rangle \bigr|
	=
	\sup_{\omegav \in \R^{\dimp} \colon \| \omegav \| = 1} 6 \bigl| \langle \Tens,\uv \otimes \omegav \otimes \omegav \rangle \bigr|
	=
	6 \| \Tens \| \, \| \uv \| .
\label{ucvkje8ij4efui2w65hvgjh}
\end{EQA}

\Subsection{Moments of a Gaussian 3-tensor}
Consider a Gaussian 3-tensor \( \Tens(\gaussv) = \langle \Tens,\gaussv^{\otimes 3} \rangle \).
Define
\begin{EQA}
	\trT_{i}
	&=&
	\sum_{j=1}^{\dimp} \Tens_{i,j,j} 
	=
	\tr \Tens_{i} \, ,
	\qquad
	i = 1,\ldots,n \, .
\label{ydy7w38uw3eiovkeredrgb}
\end{EQA}

\begin{lemma}
\label{LtensFr}
Let \( \Tens = (\Tens_{i,j,k}) \) be a 3-dimensional symmetric tensor in \( \R^{\dimp} \) and 
\( \Tens(\gaussv) = \langle \Tens, \gaussv^{\otimes 3} \rangle \) for \( \gaussv \sim \ND(0,\Id_{\dimp}) \).
With \( \trTv = (\trT_{i}) \in \R^{\dimp} \) and \( \| \Tens \|_{\Fr}^{2} = \sum_{i,j,k = 1}^{\dimp} \Tens_{i,j,k}^{2} \), 
it holds
\begin{EQA}[rcl]
	\E \bigl( \Tens(\gaussv) - 3 \langle \trTv,\gaussv \rangle \bigr)^{2}
	& = &
	6 \| \Tens \|_{\Fr}^{2} \, ,
	\\
	\E \, \Tens^{2}(\gaussv)
	&=&
	6 \| \Tens \|_{\Fr}^{2} + 9 \| \trTv \|^{2} \, .
\label{jhdyuy3e7fghy3r542tdght}
\end{EQA}
\end{lemma}

\begin{proof}
By definition
\begin{EQA}
	&& 
	\Tens(\gaussv) - 3 \langle \trTv,\gaussv \rangle 
	=
	\sum_{i,j,k=1}^{\dimp} \Tens_{i,j,k} \, \gauss_{i} \gauss_{j} \gauss_{k} 
		- 3 \sum_{i=1}^{\dimp} \gauss_{i} \, \sum_{j,k=1}^{\dimp} \Tens_{i,j,k} \, \delta_{j,k} \, .
\label{hfd7j23dmc6hwfgbo}
\end{EQA}
It is easy to see that for each \( i \) by symmetricity of \( \Tens \)
\begin{EQA}
	\E \biggl( \gauss_{i} \sum_{i,j,k=1}^{\dimp} \Tens_{i,j,k} \, \gauss_{i} \gauss_{j} \gauss_{k} \biggr)
	&=&
	3 \sum_{j \in \IIm_{i}} \Tens_{i,j,j} \E (\gauss_{i}^{2} \gauss_{j}^{2}) + \sum_{i=1}^{\dimp} \Tens_{i,i,i} \E \gauss_{i}^{4}
	=
	3 \sum_{j=1}^{\dimp} \Tens_{i,j,j} 
	=
	3 \trT_{i} \, ,
\label{vyehr6bgy4hrtdb34gybn}
\end{EQA}
where the index set \( \IIm_{i} = \{ 1,\ldots,i-1,i+1,\ldots,\dimp \} \) is obtained
by removing the index \( i \) from \( 1,\ldots,\dimp \).
This implies orthogonality 
\begin{EQA}
	\E \bigl\{ 
		\bigl( \Tens(\gaussv) - 3 \langle \trTv,\gaussv \rangle \bigr) \langle \trTv,\gaussv \rangle 
	\bigr\}
	&=&
	0 .
\label{bf6ejh7ehdsnvfetg}
\end{EQA}
Introduce the index set \( \II = \{ (i,j,k) \colon i \neq j \neq k \} \):
\begin{EQA}
	\II
	& \eqdef &
	\{ (i,j,k) \colon \Ind(i=j) + \Ind(i=k) + \Ind(j=k) = 0 \} \, .
\label{iviu38rtfgugh7rdhue}
\end{EQA}
Represent \eqref{hfd7j23dmc6hwfgbo} as
\begin{EQA}
	\Tens(\gaussv) - 3 \langle \trTv,\gaussv \rangle 
	=
	\sum_{\II} \Tens_{i,j,k} \, \gauss_{i} \, \gauss_{j} \,\gauss_{k} 
	+ 3 \sum_{i=1}^{\dimp} \sum_{j \in \IIm_{i}} \Tens_{i,j,j} \, \gauss_{i} (\gauss_{j}^{2} - 1) 
	+ \sum_{i=1}^{\dimp} \Tens_{i,i,i} \, (\gauss_{i}^{3} - 3 \gauss_{i}) \, .
\label{ij3bhj47nedtyvyhwed}
\end{EQA}
All terms in the right hand-side are orthogonal to each other allowing to compute 
\( \E \bigl( \Tens(\gaussv) - 3 \langle \trTv,\gaussv \rangle \bigr)^{2} \):
\begin{EQA}
	&& \nquad\nquad
	\E \bigl( \Tens(\gaussv) - 3 \langle \trTv,\gaussv \rangle \bigr)^{2}
	=
	\E \biggl( \sum_{\II} \Tens_{i,j,k} \, \gauss_{i} \, \gauss_{j} \,\gauss_{k} \biggr)^{2}
	\\
	&&
	+ \, \E \biggl( 
		3 \sum_{i=1}^{\dimp} \sum_{j \in \IIm_{i}} \Tens_{i,j,j} \, \gauss_{i} (\gauss_{j}^{2} - 1) 
	\biggr)^{2}
	+ \E \biggl( \sum_{i=1}^{\dimp} \Tens_{i,i,i} \, (\gauss_{i}^{3} - 3 \gauss_{i}) \biggr)^{2} .
\label{ij3bhjvu37eg47ned2wed}
\end{EQA}
Further, by symmetricity of \( \Tens \)
\begin{EQA}
	&& \nquad
	\E \biggl( \sum_{\II} \Tens_{i,j,k} \, \gauss_{i} \, \gauss_{j} \,\gauss_{k} \biggr)^{2}
	=
	\E \biggl( \sum_{\II} \Tens_{i,j,k} \, \gauss_{i} \, \gauss_{j} \,\gauss_{k} 
	\sum_{\II} \Tens_{\ic,\jc,\kc} \, \gauss_{\ic} \, \gauss_{\jc} \, \gauss_{\kc} \biggr) 
	\\
	&=&
	\E \biggl( \sum_{\II} \Tens_{i,j,k} \, \gauss_{i} \, \gauss_{j} \,\gauss_{k} 
	\sum_{(\ic,\jc,\kc) = \perm(i,j,k)} \Tens_{\ic,\jc,\kc} \, \gauss_{\ic} \, \gauss_{\jc} \, \gauss_{\kc} \biggr)
	=
	6 \sum_{\II} \Tens_{i,j,k}^{2} \, .
\label{hvhfy63trgy6rhdfybhe}
\end{EQA}
Similarly
\begin{EQA}
	\E \biggl( 3 \sum_{i=1}^{\dimp} \sum_{j \in \IIm_{i}} \Tens_{i,j,j} \, \gauss_{i} (\gauss_{j}^{2} - 1) \biggr)^{2}
	&=&
	9 \sum_{i=1}^{\dimp} \sum_{j \in \IIm_{i}} \Tens_{i,j,j}^{2} \, \E \bigl\{ \gauss_{i}^{2} (\gauss_{j}^{2} - 1)^{2} \bigr\}
	=
	18 \sum_{i=1}^{\dimp} \sum_{j \in \IIm_{i}} \Tens_{i,j,j}^{2} \, ,
	\\
	\E \biggl( \sum_{i=1}^{\dimp} \Tens_{i,i,i} \, (\gauss_{i}^{3} - 3 \gauss_{i}) \biggr)^{2}
	&=&
	\sum_{i=1}^{\dimp} \Tens_{i,i,i}^{2} \, \E (\gauss_{i}^{3} - 3 \gauss_{i})^{2}
	=
	6 \sum_{i=1}^{\dimp} \Tens_{i,i,i}^{2}
\label{h6heiubhrhdybferyfn}
\end{EQA}
yielding again by symmetricity of \( \Tens \)
\begin{EQA}
	\E \bigl( \Tens(\gaussv) - 3 \langle \trTv,\gaussv \rangle \bigr)^{2}
	&=&
	6 \sum_{\II} \Tens_{i,j,k}^{2}  
	+ 18 \sum_{i=1}^{\dimp} \sum_{j \in \IIm_{i}} \Tens_{i,j,j}^{2} 
	+ 6 \sum_{i=1}^{\dimp} \Tens_{i,i,i}^{2}
	=
	6 \| \Tens \|_{\Fr}^{2}
\label{jche6e6ghtur5ewgby32}
\end{EQA}
and assertion \eqref{jhdyuy3e7fghy3r542tdght} follows in view of orthogonality \eqref{bf6ejh7ehdsnvfetg}.
\end{proof}

Similarly we study the moments of the scaled gradient vector 
\begin{EQA}
	\TensG
	&=&
	\frac{1}{3} \nabla \Tens(\gaussv) .
\label{gfe98bm3y67vnjnh90t}
\end{EQA}
The entries \( \TensG_{i} \) of \( \TensG \) can be written as 
\( \TensG_{i} = \gaussv^{\T} \Tens_{i} \, \gaussv \); see \eqref{ho5r9hkrh4ygj32fjur}.

\begin{lemma}
\label{LmomGtens}
It holds \( \E \TensG = \trTv \),
\begin{EQA}
	\Var(\TensG)
	&=&
	\VBH^{2}
	=
	\bigl( 2 \langle \Tens_{i},\Tens_{\ic} \rangle \bigr)_{i,\ic=1,\ldots,\dimp} \, ,
\label{u7fyv65eyte56vt6rjikoy}
	\\
	\tr \VBH^{2}
	&=&
	2 \sum_{i=1}^{\dimp} \| \Tens_{i} \|_{\Fr}^{2}
	=
	2 \sum_{i,j,k=1}^{\dimp} \Tens_{i,j,k}^{2}
	=
	2 \| \Tens \|_{\Fr}^{2} \, ,
\label{dyuwe7cvjw24r5cvtjw2b21tc}
	\\
	\E \| \TensG \|^{2}
	&=&
	\| \trTv \|^{2} + 2 \| \Tens \|_{\Fr}^{2} 
	\leq 
	\frac{1}{3} \E \Tens^{2}(\gaussv) \, .
\label{udc6663yv7ryefywjwe7}
\end{EQA}
Moreover, for any \( \uv \in \R^{\dimp} \)
\begin{EQA}
	\| \VBH \uv \|^{2}
	&=&
	2 \bigl\| \Tens[\uv] \bigr\|_{\Fr}^{2} \, .
	\qquad
\label{yvgk73dhqe12q5rf64eyu}
\end{EQA}
\end{lemma}

\begin{proof}
The first statement follows directly from 
\( \E \TensG_{i} = \E \gaussv^{\T} \Tens_{i} \, \gaussv = \tr \Tens_{i} \).
For any \( i,\ic \), it holds in view of  \( \E (\gauss_{j} \gauss_{k} - \delta_{j,k})^{2} = 1 + \delta_{j,k} \) 
for all \( j,k \leq \dimp \)
\begin{EQA}
	\E (\TensG_{i} - \E \TensG_{i}) (\TensG_{\ic} - \E \TensG_{\ic})
	&=&
	\E \left( 
		\sum_{j,k=1}^{\dimp} \Tens_{i,j,k} \, (\gauss_{j} \gauss_{k} - \delta_{j,k}) 
		\sum_{\jc,\kc=1}^{\dimp} \Tens_{\ic,\jc,\kc} \, (\gauss_{\jc} \gauss_{\kc} - \delta_{\jc,\kc}) 
	\right)
	\\
	&=&
	2 \sum_{j,k=1}^{\dimp} \Tens_{i,j,k} \, \Tens_{\ic,j,k} 
	=
	2 \langle \Tens_{i},\Tens_{\ic} \rangle .
\label{ucviw38vbj4r78fu2122}
\end{EQA}
This yields \eqref{u7fyv65eyte56vt6rjikoy}.
Further
\begin{EQA}
	\tr \VBH^{2}
	&=&
	2 \sum_{i=1}^{\dimp} \langle \Tens_{i},\Tens_{i} \rangle
	=
	2 \sum_{i=1}^{\dimp} \| \Tens_{i} \|_{\Fr}^{2}
	\eqdef
	2 \| \Tens \|_{\Fr}^{2} \, ,
\label{kjv8r54jt7e367e76fvhuffg}
\end{EQA}
which proves \eqref{dyuwe7cvjw24r5cvtjw2b21tc}.
Similarly, for any \( \uv = (u_{i}) \in \R^{\dimp} \)
\begin{EQA}
	\| \VBH \uv \|^{2}
	&=&
	\uv^{\T} \VBH^{2} \uv
	=
	2 \sum_{i,\ic=1}^{\dimp} u_{i} \, u_{\ic} \langle \Tens_{i},\Tens_{\ic} \rangle
	=
	2 \biggl\| \sum_{i=1}^{\dimp} u_{i} \Tens_{i} \biggr\|_{\Fr}^{2}
	=
	2 \bigl\| \Tens[\uv] \bigr\|_{\Fr}^{2} \, 
	\qquad
\label{yvgk73dhqe12q5rf64eyupr}
\end{EQA}
completing the proof.
\end{proof}

\Subsection{\( \ell_{3}-\ell_{2} \) condition}
This section introduces a special \( \ell_{3}-\ell_{2} \) condition for a symmetric 3-tensor \( \Tens \).

\begin{description}
    \item[\label{l2l3Tref} \( \bb{(\TG)} \)]
      \textit{For some symmetric \( \dimp \)-matrix \( \TG \) and \( \tensco > 0 \),
      \( \Tens(\uv) = \langle \Tens, \uv^{\otimes 3} \rangle \) fulfills}
\begin{EQA}
	|\Tens(\uv)|
	& \leq &
	\tensco \, \| \TG \uv \|^{3} ,
	\qquad
	\uv \in \R^{\dimp} \, .
\label{7cmvvc7e3hghjj856uiedT}
\end{EQA}
\end{description}


\begin{lemma}
\label{LTensTGm}
Suppose that the tensor \( \Tens \) satisfies \nameref{l2l3Tref}.
Then
\begin{EQA}
	|\langle \Tens,\uv_{1} \otimes \uv_{2} \otimes \uv_{3} \rangle|
	& \leq &
	\tensco \, \| \TG \uv_{1} \| \, \| \TG \uv_{2} \| \, \| \TG \uv_{3} \| \, ,
	\quad
	\uv_{1} \, , \uv_{2} \, , \uv_{3} \in \R^{\dimp} \, ,
	\qquad
\label{h8wxuew7fje4uyruej}
\end{EQA}
and it holds for 
any \( \uv \in \R^{\dimp} \)
\begin{EQA}
\label{8cvmd3r6gfy33gj2wTG2}
	\| \nabla \Tens(\uv) \|
	& \leq &
	3 \tensco \, \| \TG \uv \|^{2} \, \| \TG \| \, ,
	\\
	\Tens[\uv]
	& \leq &
	\tensco \, \| \TG \uv \| \, \TG^{2} \, ,
\label{8cvmd3r6gfy33gj2wfr}
\end{EQA}
yielding
\begin{EQ}[rcl]
	\| \Tens[\uv] \|_{\Fr}^{2}
	& \leq &
	\tensco^{2} \, \| \TG \uv \|^{2} \, \tr (\TG^{4}) \, ,
	\qquad
	\uv \in \R^{\dimp} \, ,
	\\
	\| \Tens \|_{\Fr}^{2}
	& \leq &
	\tensco^{2} \, \tr (\TG^{2}) \, \tr (\TG^{4}) \, .
\label{bhdctsweghfnjggye2qv}
\end{EQ}
Further, for \( \trTv = (\trT_{i}) \in \R^{\dimp} \) with \( \trT_{i} = \tr \Tens_{i} \), it holds
\begin{EQA}
	\| \trTv \|
	& \leq &
	\tensco \, \| \TG \| \, \tr (\TG^{2}) \, ,
\label{hsfd87ewyiewy7er84eewe}
\end{EQA}
The matrix \( \VBH^{2} \) from \eqref{u7fyv65eyte56vt6rjikoy} fulfills
\begin{EQA}
	\VBH^{2}
	& \leq &
	2 \tensco^{2} \tr (\TG^{4}) \, \TG^{2} .
\label{yv8deub78546thi594}
\end{EQA}
It holds for the Gaussian tensor \( \Tens(\gaussv) \) 
\begin{EQA}
	\E \, \Tens^{2}(\gaussv)
	& \leq &
	6 \tensco^{2} \, \tr (\TG^{2}) \, \tr (\TG^{4}) + 9 \tensco^{2} \, \| \TG \|^{2} \, \tr^{2} (\TG^{2})
	\leq 
	15 \tensco^{2} \, \| \TG \|^{2} \, \tr^{2} (\TG^{2}) .
	\qquad
\label{uvjkv62hvub883jdt2jbi}
\end{EQA}
\end{lemma}

\begin{proof}
Define 3-tensor \( \Tens_{\TG} \) by \( \Tens_{\TG}(\uv) = \Tens(\TG^{-1} \uv) \). 
Then condition \eqref{7cmvvc7e3hghjj856uiedT} reads \( |\Tens_{\TG}(\uv)| \leq \tensco \) for all \( \| \uv \| \leq 1 \) while
\eqref{h8wxuew7fje4uyruej} can be written as
\begin{EQA}
	|\langle \Tens_{\TG},\uv_{1} \otimes \uv_{2} \otimes \uv_{3} \rangle|
	& \leq &
	\tensco \, ,
	\quad
	\forall \| \uv_{j} \| \leq 1, \, \, \,
	j=1,2,3.
\label{67che4v5egv5yrr9gjsfyte}
\end{EQA}
The latter holds by Banach's characterization as in \eqref{iv7dcu3921kfvhw7fjew}. 
Further, 
\begin{EQA}
	\| \nabla \Tens(\uv) \|
	&=&
	\sup_{\| \uv_{1} \| = 1} \bigl| \langle \nabla \Tens(\uv),\uv_{1} \rangle \bigr| 
	= 
	\sup_{\| \uv_{1} \| = 1} 3 \bigl| \langle \Tens,\uv \otimes \uv \otimes \uv_{1} \rangle \bigr|
	\\
	& \leq &
	3 \tensco \, \| \TG \uv \|^{2} \, \sup_{\| \uv_{1} \| = 1} \| \TG \uv_{1} \|
	\leq 
	3 \tensco \, \| \TG \uv \|^{2} \, \| \TG \| \, ,
	\\
	\| \Tens[\uv] \|
	&=&
	\sup_{\| \uv_{1} \| = 1} \langle \Tens[\uv],\uv_{1}^{\otimes 2} \rangle 
	= 
	\sup_{\| \uv_{1} \| = 1} \langle \Tens,\uv \otimes \uv_{1} \otimes \uv_{1} \rangle
	\\
	& \leq &
	\tensco \, \| \TG \uv \| \, \sup_{\| \uv_{1} \| = 1} \| \TG \uv_{1} \|^{2}
	\leq 
	\tensco \, \| \TG \uv \| \, \| \TG^{2} \| \, ,
\label{v87buib73jhf7yvbguytgu8}
\end{EQA} 
yielding \eqref{8cvmd3r6gfy33gj2wfr}.
%
Further, \( \langle \trTv, \uv \rangle = \tr \Tens[\uv] \) and by \eqref{8cvmd3r6gfy33gj2wfr}
\begin{EQA}
	\| \trTv \|
	&=&
	\sup_{\| \uv \| = 1} \bigl| \langle \trTv,\uv \rangle \bigr|
	=
	\sup_{\| \uv \| = 1} \bigl| \tr \Tens[\uv] \bigr| 
	\leq 
	\tensco \, \| \TG \| \, \tr (\TG^{2}) \, .
\label{hsfd87ewyiewy7er84eewep}
\end{EQA}
Similarly for \( \uv \in \R^{\dimp} \)
\begin{EQA}
	\bigl\| \Tens [\uv] \bigr\|_{\Fr}^{2}
	&=&	
	\tr (\Tens[\uv]^{2})
	\leq 
	\tensco^{2} \, \| \TG \uv \|^{2} \, \tr (\TG^{4}) \, .
\label{bhdctsweghfnjggye2qvp}
\end{EQA}
Finally, the use of \( \Tens_{i} = \Tens[\ev_{i}] \) for the canonic basis vectors \( \ev_{i} \in \R^{\dimp} \) yields
\begin{EQA}
	\| \Tens \|_{\Fr}^{2}
	&=&
	\sum_{i=1}^{\dimp} \tr (\Tens[\ev_{i}]^{2})
	\leq 
	\tensco^{2} \, \sum_{i=1}^{\dimp} \| \TG \ev_{i} \|^{2} \, \tr (\TG^{4})
	=
	\tensco^{2} \, \tr (\TG^{2}) \, \tr (\TG^{4}) \, ,
\label{bhdctsweghfnjggye2qvp}
\end{EQA}
and \eqref{bhdctsweghfnjggye2qv} follows.
By \eqref{yvgk73dhqe12q5rf64eyu} and \eqref{8cvmd3r6gfy33gj2wfr}, 
it holds for any \( \uv \in \R^{\dimp} \)
\begin{EQA}
	\| \VBH \uv \|^{2}
	&=&
	2 \bigl\| \Tens[\uv] \bigr\|_{\Fr}^{2}
	\leq 
	2 \tensco^{2} \tr (\TG^{4}) \| \TG \uv \|^{2} \, .
\label{yfjf7b6rrtgubuhtweee2w}
\end{EQA}
This yields \eqref{yv8deub78546thi594}.
The obtained bounds lead to \eqref{uvjkv62hvub883jdt2jbi} in view of \eqref{jhdyuy3e7fghy3r542tdght}.
\end{proof}

\Subsection{Colored case}
\label{Stenscolor}
This section extends the established upper bound to the case when the standard Gaussian vector \( \gaussv \)
is replaced by a general zero mean Gaussian vector \( \GaussD \sim \ND(0,\DPTG^{-2}) \) for 
a symmetric covariance matrix \( \DPTG^{2} \).
Then \( \GaussD = \DPTG^{-1} \gaussv \) with \( \gaussv \) standard normal and 
\( \Tens(\GaussD) = \Tens(\DPTG^{-1} \gaussv) = \Tenst(\gaussv) \) with \( \Tenst(\uv) = \Tens(\DPTG^{-1} \uv) \).
If \( \Tens \) satisfies \nameref{l2l3Tref} then \( \Tenst \) does as well
but \( \TG^{2} \) has to be replaced by \( \TGD^{2} = \DPTG^{-1} \TG^{2} \DPTG^{-1} \).

\begin{lemma}
\label{LTensTGmG}
Let \( \Tens(\uv) \) satisfies \nameref{l2l3Tref} with some \( \TG \) and \( \tensco \).
Then \( \Tenst(\uv) = \Tens(\DPTG^{-1} \uv) \) satisfies \nameref{l2l3Tref} with 
\( \TGD^{2} = \DPTG^{-1} \TG^{2} \DPTG^{-1} \) in place of \( \TG^{2} \) and the same \( \tensco \).
In particular, with \( \Tenst = (\Tenst_{i}) \), \( \trTvt = (\tr \Tenst_{i}) \), and 
\( \VBHt^{2} \eqdef \bigl( 2 \langle \Tenst_{i},\Tenst_{\ic} \rangle \bigr)_{i,\ic=1,\ldots,\dimp} \), it holds
\begin{EQA}
\label{bhdctsweghfnjggye2qvt}
	\| \Tenst \|_{\Fr}^{2}
	& \leq &
	\tensco^{2} \, \tr (\TGD^{2}) \, \tr (\TGD^{4}),
	\\
	\| \trTvt \|
	& \leq &
	\tensco \, \| \TGD \| \, \tr (\TGD^{2}) \, ,
\label{7vjnvc7e3yfgb6ywvbueom}
	\\
	\VBHt^{2}
	& \leq &
	2 \tensco^{2} \, \tr (\TGD^{4}) \, \TGD^{2} ,
\label{7vjnvc7e3yfgb6ywvbueuj}
\end{EQA}
Moreover, for any \( \uv \in \R^{\dimp} \) 
\begin{EQ}[rcl]
	\| \nabla \Tenst(\uv) \|
	& \leq &
	3 \tensco \, \| \TGD \, \uv \|^{2} \, \| \TGD \| \, .
	\\
	\bigl\| \Tenst[\uv] \bigr\|_{\Fr}^{2}
	& \leq &
	\tensco^{2} \, \| \TGD \, \uv \|^{2} \, \tr \TGD^{4} \, ,
\label{vdf8we3hjg7gt4f7yvbg}
\end{EQ}
\end{lemma}

\begin{proof}
By definition, for any \( \uv \in \R^{\dimp} \)
\begin{EQA}
	\Tenst(\uv)
	&=&
	\Tens(\DPTG^{-1} \uv)
	\leq 
	\tensco \| \TG \DPTG^{-1} \uv \|^{3} 
	=
	\tensco \| \TGD \, \uv \|^{3} 
\label{d7vuvuufu7e4e347nsduyc}
\end{EQA}
yielding \nameref{l2l3Tref} for \( \Tenst \).
Now Lemma~\ref{LTensTGmG} enables us to apply the results of 
Lemma~\ref{LTensTGm} with \( \TGD \) in place of \( \TG \).
Finally, for any \( \uv \) with \( \| \TGD \, \uv \| \leq \rrTG \), it holds by \eqref{8cvmd3r6gfy33gj2wTG2}
\begin{EQA}
	\| \nabla \Tenst(\uv) \|
	& \leq &
	3 \tensco \, \| \TGD \, \uv \|^{2} \, \| \TGD \|
	\leq 
	\tensco \, \rrTG^{2} \, \| \TGD \| \, .
\label{f8uvuvuv7e7edgw2hgc7ww}
\end{EQA}
This completes the proof.
\end{proof}

\begin{lemma}
\label{LcolorGauss}
For the Gaussian tensor \( \Tens(\GaussD) \), it holds 
\begin{EQA}
	\E \, \Tens^{2}(\GaussD)
	& \leq &
	15 \tensco^{2} \, \| \TGD \|^{2} \, \tr^{2} (\TGD^{2}) .
\label{uvjkb7nh3f6eyeyudt2jbi}
\end{EQA}
Moreover, with some fixed value \( \const_{2} \), it holds
\begin{EQA}
	\sqrt{\E \, \Tens^{4}(\GaussD)}
	& \leq &
	\const_{2} \, \tensco^{2} \, \| \TGD \|^{2} \, \tr^{2} (\TGD^{2}) .
\label{uvjkb7nh3f6eyeyudt2jbi4}
\end{EQA}
\end{lemma}

\begin{proof}
Apply Lemma~\ref{LtensFr} to \( \Tens(\GaussD) = \Tenst(\gaussv) \) 
and use \eqref{bhdctsweghfnjggye2qvt} and \eqref{7vjnvc7e3yfgb6ywvbueom}.
The bound for \( \E \, \Tens^{4}(\GaussD) \) can be derived from Lemma~\ref{Lmomentsexp} later
with an explicit constant \( \const_{2} \). 
\end{proof}

\Subsection{An exponential bound on \( \Tens(\GaussD) \)}
Let \( \dlt(\uv) \) be a smooth function on \( \R^{\dimp} \).
A typical example we have in mind is \( \dlt(\uv) = \Tens(\uv) \), where
\( \Tens \) is a symmetric 3-tensor satisfying \nameref{l2l3Tref} with some \( \TG^{2} \) and \( \tensco \).
Let also \( \GaussD \sim \ND(0,\DPTG^{-2}) \).
Our aim is a possibly accurate exponential bound for \( \dlt(\GaussD) \), in particular, for the Gaussian tensor  
\( \Tens(\GaussD) = \langle \Tens,\GaussD^{\otimes 3} \rangle \).
We use \( \dlt(\GaussD) = \dlt(\DPTG^{-1} \gaussv) = \dltt(\gaussv) \) for \( \dltt(\uv) = \dlt(\DPTG^{-1} \uv) \)
and \( \gaussv \) standard normal.
The results use a bound on the norm \( \| \nabla \dltt(\uv) \| \) which is hard to verify on the whole domain \( \R^{\dimp} \).
Therefore, we limit the domain of \( \dlt(\uv) \) to a subset \( \UV \) on which the Gaussian measure 
\( \ND(0,\DPTG^{-2}) \) well concentrates.
Clearly, for any \( \uv \in \R^{\dimp} \) 
\begin{EQA}
	\nabla \dltt(\uv) 
	&=&
	\DPTG^{-1} \nabla \dlt(\DPTG^{-1} \uv) .
\label{uv7yrjrf846yr9df}
\end{EQA}
For the rest of this section, we assume 
that \( \| \nabla \dltt(\uv) \| \) is uniformly bounded over the set of \( \uv \) with \( \DPTG^{-1} \uv \in \UV \)
for the set \( \UV \) with the elliptic shape 
\begin{EQA}
	\UV
	& \eqdef &
	\{ \uv \colon \| \TG \uv \| \leq \rr \} .
\label{fc7vhfe367hf78r78ewhc}
\end{EQA}
This applies to \( \dlt(\uv) = \Tens(\uv) \) for a tensor \( \Tens \) satisfying \nameref{l2l3Tref}.
We consider local behavior of \( \dlt(\GaussD) \) for \( \GaussD \sim \ND(0,\DPTG^{-2}) \).
With \( \UV \) fixed, introduce the notation 
\begin{EQA}[c]
	\EUV \, \xi \eqdef \E \{ \xi \Ind(\GaussD \in \UV) \} .
\label{yvjby6wehybnbuejwfct}
\end{EQA}
Remind the definition \( \TGD^{2} =	\DPTG^{-1} \TG^{2} \DPTG^{-1} \).
The next lemma explains the choice of the radius \( \rr \) to ensure a concentration effect of \( \GaussD \) on \( \UV \).

\begin{lemma}
\label{PTensG}
For a fixed \( \xx \), set \( \rr = \rr(\xx) = \zq(\TGD^{2},\xx) \) with   
\begin{EQA}
	\zq^{2}(\TGD^{2},\xx) 
	&=& 
	\tr(\TGD^{2}) + 2 \sqrt{\xx \tr(\TGD^{4})} + 2 \xx \| \TGD^{2} \| .
\label{fubnv6e67ehfvyterbue3u}
\end{EQA}
For the set \( \UV \) from \eqref{fc7vhfe367hf78r78ewhc}, suppose
\begin{EQA}
	\sup_{\vv \colon \DPTG^{-1} \vv \in \UV}
	\| \nabla \dltt(\vv) \|
	=
	\sup_{\uv \in \UV}
	\| \DPTG^{-1} \nabla \dlt(\uv) \|
	& \leq &
	\grad \, .
\label{yfuejfivi8eujnbwerwree}
\end{EQA}
Then it holds for \( X = \dlt(\GaussD) - \EUV \, \dlt(\GaussD) \), 
with any \( \muH \) and any integer \( k \)
\begin{EQA}[ccll]
	\EUV \, \ex^{\muH X} 
	& \leq &
	\exp\bigl( \muH^{2} \grad^{2}/2 \bigr) \, ,
	&
\label{v7jdciuwmedfib5i84k}
	\\
	\EUV |X|^{2k} 
	& \leq &
	\CONSTi_{k}^{2} \grad^{2k} \, ,
	&
	\CONSTi_{k}^{2} = 2^{k+1} k! \, .
\label{yjwvue4ufvnwyvne}
\end{EQA}
Also
\begin{EQA}
	\P\Bigl( X > \grad \sqrt{2\xx} \Bigr)
	& \leq &
	2 \ex^{-\xx} .
\label{76djhv83iurfgjhyhP}
\end{EQA}
\end{lemma}

\begin{proof}
With \( \gaussv \sim \ND(0,\Id_{\dimp}) \) and \( \GaussD \sim \ND(0,\DPTG^{-2}) \), it holds
\begin{EQA}
	\P(\GaussD \not\in \UV)
	&=&
	\P\bigl( \| \TG \DPTG^{-1} \gaussv \| > \rr \bigr) .
\label{tduc7wquyjhiw6wyhdj}
\end{EQA}
For \( \rr = \zq(\TGD^{2},\xx) \), 
Gaussian concentration bound yields
\begin{EQA}
	\P(\GaussD \not\in \UV)
	=
	\P\bigl( \| \TG \DPTG^{-1} \gaussv \| > \zq(\TGD^{2},\xx) \bigr)
	& \leq &
	\ex^{-\xx}.
\label{xc7y6wqjhc9iwndreae}
\end{EQA}
Further, \( \dlt(\GaussD) = \dlt(\DPTG^{-1} \gaussv) = \dltt(\gaussv) \) for \( \gaussv \) standard normal
and by \eqref{yfuejfivi8eujnbwerwree}, the norm of the gradient \( \nabla \dltt(\vv) \) is bounded by \( \grad \) 
for all \( \vv \) with \( \DPTG^{-1} \vv \in \UV \).
The use of log-Sobolev inequality and Herbst's arguments yields \eqref{v7jdciuwmedfib5i84k} for 
\( X = \dlt(\DPTG^{-1} \gaussv) - \EUV \dlt(\DPTG^{-1} \gaussv) \);
see Theorem 5.5 in \cite{boucheron2013concentration} or Proposition 5.4.1 in \cite{bakry2013analysis}.
Result \eqref{v7jdciuwmedfib5i84k} also implies the probability bound 
\begin{EQA}
	\P\Bigl( X > \sqrt{2\xx} \, \grad \Bigr)
	& \leq &
	\P(\GaussD \not\in \UV) + \P\Bigl( X > \grad \sqrt{2\xx} \, , \, \GaussD \in \UV \Bigr)
	\leq 
	2 \ex^{-\xx} ;
\label{76djhv83iurfgjhyhy4r}
\end{EQA}
see (5.4.2) in \cite{bakry2013analysis}.
Now Lemma~\ref{Lmomentsexp} and \eqref{v7jdciuwmedfib5i84k} imply \eqref{yjwvue4ufvnwyvne}.
\end{proof}

\begin{lemma}[\cite{boucheron2013concentration}, Theorem 2.1]
\label{Lmomentsexp}
Let a r.v. \( X \) satisfy 
\( \E \exp (\muH X) \leq \exp(\muH^{2} \grad^{2}/2) \) for all \( \muH \) with some \( \grad^{2} > 0 \).
Then for any integer \( k \)
\begin{EQA}
	\E |X|^{2k}
	& \leq &
	\CONSTi_{k}^{2} \, \grad^{2k} \, ,
	\qquad
	\CONSTi_{k}^{2} = 2^{k+1} k! 
	\, .
\label{6cyvhwyvhejhewycnwsbj}
\end{EQA}
In particular, \( \CONSTi_{1} = 2 \), \( \CONSTi_{2} = 4 \), \( \CONSTi_{3} = \sqrt{96} \leq 10 \),
\( \CONSTi_{4} = 16 \sqrt{3} \leq 28 \).
\end{lemma}

\begin{proof}
Conditions of the lemma and Markov inequality imply for any \( u > 0 \) with \( \mu = u \)
\begin{EQA}
	\P\bigl( X/ \grad > u \bigr)
	& \leq &
	\ex^{- \mu u} \E \exp (\muH X/\grad)
	\leq 
	\exp(- u^{2} /2)
\label{gsnvfu3ruhadsqsgvdfw}
\end{EQA}
and similarly for \( \P\bigl( - X/ \grad > u \bigr) \) hence, 
\begin{EQA}
	\E |X/s|^{2k} 
	&=&
	\int_{0}^{\infty} \P\bigl( |X/s|^{2k} > x \bigr) \, dx
	=
	2k \int_{0}^{\infty} x^{2k-1} \P\bigl( |X/s| > x \bigr) \, dx
	\\
	& \leq &
	4k \int_{0}^{\infty} x^{2k-1} \ex^{-x^{2}/2} \, dx
	= 
	4k \int_{0}^{\infty} (2t)^{k-1} \ex^{-t} \, dt
	=
	2^{k+1} k! 
\label{jvt6wehbf7y67e34nhfdd}
\end{EQA}
as claimed. 
\end{proof}

Let \( X \) satisfy \( \E \exp (\muH X) \leq \exp(\muH^{2} \grad^{2}/2) \) for all \( \muH \) with some \( \grad \) small.
One can expect that \( \ex^{X} \) can be well approximated for \( k \geq 2 \) by
\begin{EQA}
	\EX_{k}(X)
	& \eqdef &
	1 + X + \ldots + \frac{X^{k-1}}{(k-1)!} \, .
\label{tyfcyhwwuv7wj2apc7ew3j}
\end{EQA}

\begin{lemma}
\label{PTensG3}
Let a random variable \( X \) satisfy 
\( \E \exp (\muH X) \leq \exp(\muH^{2} \grad^{2}/2) \) for all \( \muH \) with some \( \grad^{2} > 0 \).
Then for a random variable \( \xi \) such that \( |\xi| \leq 1 \) and any integer \( k \) with \( \CONSTi_{k} \) from \eqref{yjwvue4ufvnwyvne} and \( \EX_{k}(X) \) from \eqref{tyfcyhwwuv7wj2apc7ew3j}
\begin{EQA}
	\bigl| 
		\E \bigl( \ex^{X} - \EX_{k}(X) \bigr) \xi
	\bigr|
	& \leq &
	\frac{\CONSTi_{k}}{k!} \, \grad^{k} \ex^{\grad^{2}} \, .
\label{7vnf53g6fy363ghfughd}
\end{EQA}
In particular, with \( \CONSTi_{2} = 4 \) and \( \CONSTi_{3}^{2} = 96 \)
\begin{EQ}[rcl]
	\bigl| 
		\E \bigl( \ex^{X} - 1 - X \bigr) \xi 
	\bigr|
	& \leq &
	2 \grad^{2} \, \ex^{\grad^{2}} \, ,
	\\
	\bigl| 
		\E \bigl( \ex^{X} - 1 - X - \frac{X^{2}}{2} \bigr) \xi 
	\bigr|
	& \leq &
	\frac{5}{3} \grad^{3} \, \ex^{\grad^{2}} \, .
\label{7vnf53g6fy363ghfughd2}
\end{EQ}
If \( \xi \) is not bounded but \( \E \xi^{2k+2} < \infty \), then with \( \rho = k/(k+1) \) 
\begin{EQA}
	\bigl| 
		\E \bigl( \ex^{X} - \EX_{k}(X) \bigr) \xi 
	\bigr|
	& \leq &
	\frac{\CONSTi_{k+1}^{\rho}}{k!} \, \grad^{k} \, \ex^{\grad^{2}} \, \bigl( \E \xi^{2k+2} \bigr)^{\frac{1}{2k+2}} \, .
\label{7vnf53g6fy363ghfugh62d}
\end{EQA}
In particular, with \( \CONSTi_{3}^{2/3} = 96^{1/3} \leq 4.6 \) and \( \CONSTi_{4}^{3/4} \leq 12 \)
\begin{EQA}
	\bigl| 
		\E \bigl( \ex^{X} - 1 - X \bigr) \xi 
	\bigr|
	& \leq &
	2.3 \, \grad^{2} \, \ex^{\grad^{2}} \, \bigl( \E \xi^{6} \bigr)^{1/6} \, ,
\label{7vnf53g6fy363ghfugh2}
	\\
	\bigl| 
		\E \bigl( \ex^{X} - 1 - X - \frac{X^{2}}{2} \bigr) \xi 
	\bigr|
	& \leq &
	2 \, \grad^{3} \, \ex^{\grad^{2}} \, \bigl( \E \xi^{8} \bigr)^{1/8} \, .
\label{7vnf53g6fy363ghfugh4}
\end{EQA}
\end{lemma}

\begin{proof}
Define 
\begin{EQA}
	\riskt(t)
	& \eqdef &
	\E \bigl\{ \bigl( \ex^{t X} - \EX_{k}(tX) \bigr) \xi \bigr\}  \, .
\label{cuvh46fghb83bbuwbdu}
\end{EQA}
Obviously \( \riskt(0) = \riskt'(0) = \ldots = \riskt^{(k-1)}(0) = 0 \).
The Taylor expansion of order \( k \) yields
\begin{EQA}
	\bigl| \riskt(1) \bigr|
	& \leq &
	\frac{1}{k!} \sup_{t \in [0,1]} |\riskt^{(k)}(t)| .
\label{hjvg7ejfi9i4nretcrwty}
\end{EQA}
Further, 
\begin{EQA}
	\riskt^{(k)}(t)
	& = &
	\E (X^{k} \,  \xi \, \ex^{t X}) \, .
\label{hfywhjfvueybwfsgh}
\end{EQA}
Consider first the case \( |\xi| \leq 1 \) a.s.
By the Cauchy-Schwarz inequality, \eqref{v7jdciuwmedfib5i84k} of Lemma~\ref{PTensG}, and \eqref{6cyvhwyvhejhewycnwsbj},
it holds  for any \( t \in [0,1] \) 
\begin{EQA}
	|\riskt^{(k)}(t)|^{2}
	& \leq &
	\E |X|^{2k} \,\, \E \, \ex^{2 t X} 
	\leq 
	\CONSTi_{k}^{2} \grad^{2k} \ex^{2 \grad^{2}}
	\, .
\label{hgdrbketrdefvbdtctgewb}
\end{EQA}
For a general \( \xi \), in a similar way, it holds with \( \rho = k/(k+1) \)
\begin{EQA}
	|\riskt^{(k)}(t)|^{2}
	& \leq &
	\E (|X|^{2k} \xi^{2}) \,\, \E \, \ex^{2 t X} 
	\leq 
	\bigl( \E |X|^{2k+2} \bigr)^{\rho} \bigl( \E \xi^{2k+2} \bigr)^{1 - \rho} \,\, \E \, \ex^{2 \grad^{2}} 
	\\
	& \leq &
	\CONSTi_{k+1}^{2\rho} \grad^{2k} \ex^{2 \grad^{2}} \, \bigl( \E \xi^{2k+2} \bigr)^{1 - \rho} \, ,
\label{hv6gdffv5egd5wghb87ruwe}
\end{EQA}
and \eqref{7vnf53g6fy363ghfugh62d} follows.
\end{proof}

This result with \( \xi = 1 \) yields an approximation \( \E \ex^{X} \approx 1 + \E X + \E X^{2}/2 \)
and with \( \xi = X \) an approximation \( \E (X \ex^{X}) \approx \E X + \E X^{2} \).

\begin{lemma}
\label{PTensGEX}
Let a random variable \( X \) satisfy 
\( \E \exp (\muH X) \leq \exp(\muH^{2} \grad^{2}/2) \) for all \( \muH \) with some \( \grad^{2} > 0 \).
Then
\begin{EQA}
	\bigl| \E \ex^{X} - 1 - \E X - \E X^{2}/2) \bigr|
	& \leq &
	2 \grad^{3} \ex^{\grad^{2}} \, ,
	\\
	\bigl| \E (X \ex^{X}) - \E X - \E X^{2}) \bigr|
	& \leq &
	5 \grad^{3} \ex^{\grad^{2}} \, .
\label{iucfjhc76ehfyrujgbge}
\end{EQA}
\end{lemma}

\begin{proof}
The first bound follows from \eqref{7vnf53g6fy363ghfugh4} with \( \xi \equiv 1 \).
Further, \eqref{6cyvhwyvhejhewycnwsbj} for \( k=3 \) implies \( \E X^{6} \leq 96 \grad^{6} \) and
\eqref{7vnf53g6fy363ghfugh2} with \( \xi = X \) yields the second bound.
\end{proof}

Now we specify the obtained bounds for two scenarios. 
Let \( \lgd \) be a function on \( \R^{\dimp} \).
First we consider a symmetric 3-tensor \( \Tens \) which can be viewed as third order derivative of \( \lgd \) 
at some point \( \xv \) and define \( X = \Tens(\GaussD) \) for \( \GaussD \sim \ND(0,\DVL^{-2}) \).

\begin{lemma}
\label{PtensGTens}
Let \( \Tens \) be a symmetric 3-tensor \( \Tens \) satisfying \nameref{l2l3Tref} and \( \GaussD \sim \ND(0,\DPTG^{-2}) \).
Consider the set \( \UV \) from Lemma~\ref{PTensG}.
Then all the statements of Lemma~\ref{PTensG} and Lemma~\ref{PTensG3} continue to apply with \( X = \Tens(\GaussD) \)
and \( \grad \eqdef 3 \, \tensco \, \rr^{2} \| \TGD \| \).
In particular, it holds for any \( \muH \) and any integer \( k \)
\begin{EQA}[ccl]
	\EUV \, \ex^{\muH \Tens(\GaussD)} 
	& \leq &
	\exp\bigl( \muH^{2} \grad^{2}/2 \bigr) ,
\label{v7jdciuwmedfib5i84kT}
	\\
	\EUV |\Tens(\GaussD)|^{2k} 
	& \leq &
	\CONSTi_{k}^{2} \grad^{2k} \, , \qquad \qquad \CONSTi_{k} = 2^{k+1} k! \, .
\label{yjwvue4ufvnwyvneT}
\end{EQA}
\end{lemma}

\begin{proof}
Condition \nameref{l2l3Tref} ensures that 
the gradient of \( \Tens(\DVL^{-1} \uv) \) is uniformly bounded by
\( \grad = 3 \, \tensco \, \rr^{2} \| \TGD \| \) on the local set \( \UV \); see Lemma~\ref{LTensTGmG}.
This enables the statements of Lemma~\ref{PTensG} and Lemma~\ref{PTensG3}.
\end{proof}

For the second scenario, let \( \dlt(\uv) \) be  the third order remainder in the Taylor expansion of 
\( \lgd(\xv + \uv) \) at a fixed point \( \xv \):
\begin{EQA}
	\dlt(\uv)
	& \eqdef &
	\lgd(\xv + \uv) - \lgd(\xv) - \langle \nabla \lgd(\xv), \uv \rangle 
	- \frac{1}{2} \langle \nabla^{2} \lgd(\xv), \uv^{\otimes 2} \rangle 
\label{ytvcjf63wehg7ehfyjjf}
\end{EQA}
and consider \( \dlt(\GaussD) \).
In this case we assume that \( \lgd \) satisfies the following condition.
 
\begin{description}
    \item[\label{l2l3gref} \( \bb{(\TG_{3})} \)]
    \textit{For some \( \TG \) and \( \dltwu_{3} > 0 \), it holds 
     for any \( \uv \in \UV = \{ \uv \colon \| \TG \uv \| \leq \rr \} \),
      }
\begin{EQA}
	|\langle \nabla^{3} \lgd(\xv + \uv), \uv_{1}^{\otimes 3} \rangle|
	& \leq &
	\dltwu_{3} \, \| \TG \uv_{1} \|^{3} ,
	\qquad
	\uv_{1} \in \R^{\dimp} .
\label{vjb7hyer5ewre5fty4fghg}
\end{EQA}
\end{description}

Banach's characterization \cite{Banach1938} yields for any \( \uv_{1} , \uv_{2}, \uv_{3} \in \R^{\dimp} \)
\begin{EQA}
	\bigl| \langle \nabla^{3} \lgd(\xv + \uv), \uv_{1} \otimes \uv_{2} \otimes \uv_{3} \rangle \bigr|
	& \leq &	 
	\dltwu_{3} \, \| \TG \uv_{1} \| \, \| \TG \uv_{2} \| \, \| \TG \uv_{3} \| \, ;
\label{jbufiv784ejf76e3n94}
\end{EQA}
see Lemma~\ref{LTensTGm}.

\begin{lemma}
\label{PtensGdlt}
Let a function \( \lgd \) satisfy \nameref{l2l3gref} and \( \GaussD \sim \ND(0,\DPTG^{-2}) \).
Define \( \TGD \) by \( \TGD^{2} = \DPTG^{-1} \TG^{2} \DPTG^{-1} \).
Then all the statements of Lemma~\ref{PTensG} and Lemma~\ref{PTensG3} continue to apply with \( X = \dlt(\GaussD) \)
for \( \dlt(\uv) \) from \eqref{ytvcjf63wehg7ehfyjjf}
and \( \grad \eqdef \dltwu_{3} \, \rr^{2} \| \TGD \|/2 \).
\end{lemma}

\begin{proof}
Define \( \dltt(\uv) = \dlt(\DPTG^{-1} \uv) \). 
Note that \( \DPTG^{-1} \uv \in \UV \) means \( \| \TG \DPTG^{-1} \uv \| = \| \TGD \uv \| \leq \rr \).
We only have to check that condition \nameref{l2l3gref} implies with \( \grad = \dltwu_{3} \, \rr^{2} \| \TGD \|/2 \)
\begin{EQA}
	\sup_{\uv \colon \| \TGD \uv \| \leq  \rr} \| \nabla \dltt(\uv) \|
	=
	\sup_{\uv \colon \| \TGD \uv \| \leq  \rr} \| \DPTG^{-1} \nabla \dlt(\uv) \|
	& \leq &
	\grad .
\label{i8foi9eoeg8j4dfgifert9i}
\end{EQA}
Indeed, the Taylor expansion at \( \uv = 0 \) yields by \( \nabla \dlt(0) = 0 \) and \( \nabla^{2} \dlt(0) = 0 \)
\begin{EQA}
	\| \DPTG^{-1} \nabla \dlt(\uv) \|
	& \leq &
	\sup_{\| \ww \| = 1} \langle \nabla \dlt(\uv), \DPTG^{-1} \wv \rangle 
	=
	\sup_{\| \ww \| = 1} \langle \nabla \dlt(\uv) - \nabla \dlt(0) - \nabla^{2} \dlt(0) \uv, \DPTG^{-1} \wv \rangle
	\\
	&=&
	\frac{1}{2} \sup_{\| \ww \| = 1} \langle \nabla^{3} \dlt(t \uv) , \uv \otimes \uv \otimes \DPTG^{-1} \wv \rangle
\label{bn9m3gf8rjdgtgi8evig}
\end{EQA}
By \eqref{jbufiv784ejf76e3n94}
\begin{EQA}
	\| \DPTG^{-1} \nabla \dlt(\uv) \|
	& \leq &
	\frac{1}{2} \sup_{\| \ww \| = 1} \langle \nabla^{3} \dlt(t \uv) , \uv \otimes \uv \otimes \DPTG^{-1} \wv \rangle
	\leq 
	\frac{\dltwu_{3}}{2} \sup_{\| \ww \| = 1} \| \TG \uv \| \, \| \TG \uv \| \, \| \TG \DPTG^{-1} \wv \|
	\\
	& \leq &
	\frac{\dltwu_{3}}{2} \, \rr^{2} \, \sup_{\| \ww \| = 1} \, \| \TGD \wv \|
	=
	\frac{\dltwu_{3}}{2} \, \rr^{2} \, \| \TGD \| \, 
\label{hgf7ehf7ru3nmghkhucdjh}
\end{EQA}
as required.
\end{proof}